\documentclass{article}

\usepackage{amssymb,latexsym,amsmath}

\usepackage[pdftex]{graphicx}

\begin{document}

\newcommand{\bfi}{\bfseries\itshape}

\makeatletter

\@addtoreset{figure}{section}

\def\thefigure{\thesection.\@arabic\c@figure}

\def\fps@figure{h, t}

\@addtoreset{table}{bsection}

\def\thetable{\thesection.\@arabic\c@table}

\def\fps@table{h, t}

\@addtoreset{equation}{section}

\def\theequation{\thesubsection.\arabic{equation}}

\makeatother

\newtheorem{thm}{Theorem}[section]

\newtheorem{prop}[thm]{Proposition}

\newtheorem{lema}[thm]{Lemma}

\newtheorem{cor}[thm]{Corollary}

\newtheorem{defi}[thm]{Definition}

\newtheorem{rk}[thm]{Remark}

\newtheorem{exempl}{Example}[section]

\newenvironment{exemplu}{\begin{exempl}  \em}{\hfill $\square$

\end{exempl}}

\newcommand{\comment}[1]{\par\noindent{\raggedright\texttt{#1}

\par\marginpar{\textsc{Comment}}}}

\newcommand{\todo}[1]{\vspace{5 mm}\par \noindent \marginpar{\textsc{ToDo}}\framebox{\begin{minipage}[c]{0.95 \textwidth}

\tt #1 \end{minipage}}\vspace{5 mm}\par}

\newcommand{\ea}{\mbox{{\bf a}}}

\newcommand{\eu}{\mbox{{\bf u}}}

\newcommand{\ueu}{\underline{\eu}}

\newcommand{\ueo}{\overline{u}}

\newcommand{\oeu}{\overline{\eu}}

\newcommand{\ew}{\mbox{{\bf w}}}

\newcommand{\ef}{\mbox{{\bf f}}}

\newcommand{\eF}{\mbox{{\bf F}}}

\newcommand{\eC}{\mbox{{\bf C}}}

\newcommand{\en}{\mbox{{\bf n}}}

\newcommand{\eT}{\mbox{{\bf T}}}

\newcommand{\eL}{\mbox{{\bf L}}}

\newcommand{\eR}{\mbox{{\bf R}}}

\newcommand{\eV}{\mbox{{\bf V}}}

\newcommand{\eU}{\mbox{{\bf U}}}

\newcommand{\ev}{\mbox{{\bf v}}}

\newcommand{\eve}{\mbox{{\bf e}}}

\newcommand{\uev}{\underline{\ev}}

\newcommand{\eY}{\mbox{{\bf Y}}}

\newcommand{\eK}{\mbox{{\bf K}}}

\newcommand{\eP}{\mbox{{\bf P}}}

\newcommand{\eS}{\mbox{{\bf S}}}

\newcommand{\eJ}{\mbox{{\bf J}}}

\newcommand{\eB}{\mbox{{\bf B}}}

\newcommand{\eH}{\mbox{{\bf H}}}

\newcommand{\leb}{\mathcal{ L}^{n}}

\newcommand{\eI}{\mathcal{ I}}

\newcommand{\eE}{\mathcal{ E}}

\newcommand{\hen}{\mathcal{H}^{n-1}}

\newcommand{\eBV}{\mbox{{\bf BV}}}

\newcommand{\eA}{\mbox{{\bf A}}}

\newcommand{\eSBV}{\mbox{{\bf SBV}}}

\newcommand{\eBD}{\mbox{{\bf BD}}}

\newcommand{\eSBD}{\mbox{{\bf SBD}}}

\newcommand{\ecs}{\mbox{{\bf X}}}

\newcommand{\eg}{\mbox{{\bf g}}}

\newcommand{\paromega}{\partial \Omega}

\newcommand{\gau}{\Gamma_{u}}

\newcommand{\gaf}{\Gamma_{f}}

\newcommand{\sig}{{\bf \sigma}}

\newcommand{\gac}{\Gamma_{\mbox{{\bf c}}}}

\newcommand{\deu}{\dot{\eu}}

\newcommand{\dueu}{\underline{\deu}}

\newcommand{\dev}{\dot{\ev}}

\newcommand{\duev}{\underline{\dev}}

\newcommand{\weak}{\stackrel{w}{\approx}}

\newcommand{\mild}{\stackrel{m}{\approx}}

\newcommand{\lrightarrow}{\stackrel{L}{\rightarrow}}

\newcommand{\rrightarrow}{\stackrel{R}{\rightarrow}}

\newcommand{\strong}{\stackrel{s}{\approx}}

\newcommand{\weakdown}{\rightharpoondown}

\newcommand{\opg}{\stackrel{\mathfrak{g}}{\cdot}}

\newcommand{\opunu}{\stackrel{1}{\cdot}}
\newcommand{\opdoi}{\stackrel{2}{\cdot}}

\newcommand{\opn}{\stackrel{\mathfrak{n}}{\cdot}}
\newcommand{\opx}{\stackrel{x}{\cdot}}

\newcommand{\tr}{\ \mbox{tr}}

\newcommand{\Ad}{\ \mbox{Ad}}

\newcommand{\ad}{\ \mbox{ad}}

\renewcommand{\contentsname}{ }

\title{Normed groupoids with dilations}

\author{Marius Buliga\footnote{"Simion Stoilow" Institute of Mathematics of the 
Romanian Academy, PO BOX 1-764,014700 Bucharest, Romania, e-mail: 
Marius.Buliga@imar.ro }}

\date{This version:  14.07.2011}

\maketitle

\begin{abstract}
We study normed groupoids with dilations and their induced deformations. 
\end{abstract}


\tableofcontents


\vspace{1.cm}

\section{Normed groupoids}

\subsection{Groupoids}

A groupoid is a small category whose arrows are all invertible. More 
precisely we have the following definition. 

\begin{defi}
A {\bf groupoid} over a set $X$ is a set of {\bf arrows} $G$ along with 
 a {\bf target map} $\omega: G \rightarrow X$, 
a {\bf source map}  $\alpha: G \rightarrow X$,  a {\bf identity section} 
$e: X \rightarrow G$ which is a injective function,  a {\bf partially defined operation} (or product)
 $m$ on $G$, which  is a function: 
 $$m: G^{(2)} \, = \, \left\{ (g,h) \in G^{2} 
\mbox{ : } \omega(h) = \alpha(g) \right\} \rightarrow G \, , \quad m(g,h) \, = \,
gh$$
and a {\bf inversion map} $\displaystyle inv: G \rightarrow G$, $\displaystyle 
inv(g) = g^{-1}$. 
These are the structure maps of the groupoid. They satisfy  several identities. 
\begin{enumerate}
\item[(a)] 
For any $\displaystyle (g,h) \in G^{(2)}$ we
have 
$$\omega(gh) \, = \, \omega(g) \, , \quad  \alpha(gh) \, = \, \alpha(h)$$ 
\item[(b)] Then 
for any $\displaystyle (g,h), (gh,k) \in G^{(2)}$  we have also 
$\displaystyle (g, hk) \in G^{(2)}$. This allows us to write the expression 
$g(hk)$ and to state that the operation $m$ is associative: 
$$ (gh)k \, = \, g(hk)$$
\item[(c)]  for any $g \in G$ the identity section satisfies 
$\displaystyle (e(\omega(g)), g),
(g, e(\alpha(g))) \in G^{(2)})$ and 
$$e(\omega(g)) \, g \, = \,  g\, e(\alpha(g)) \, = \, g$$
\item[(c)] The inversion map is an involution: $inv \, inv \, = \, id$. For any 
$g \in G$ we have $\displaystyle (g^{-1}, g)\in G^{(2)}$ and 
$\displaystyle (g , g^{-1}) \in G^{(2)}$ and 
$$g^{-1} \, g \, = \, e(\alpha(g)) \, , \quad g \, g^{-1} \, = \, e(\omega(g))$$
\end{enumerate}
\label{defgroid}
\end{defi}

An equivalent definition of a groupoid emphasizes the fact that a groupoid may
be defined only in terms of its arrows. 

\begin{defi}
A groupoid is a set $G$ with two operations $\displaystyle inv: G \rightarrow 
G$,  $\displaystyle m: G^{(2)} \subset G \times
G \rightarrow G$, which satisfy a number of properties. With the notations 
$\displaystyle inv(a) = a^{-1}$, $\displaystyle m(a,b) = ab$, these properties
are: for any $a,b,c \in G$ 
\begin{enumerate}
\item[(i)] if $\displaystyle (a,b) \in G^{(2)}$ and 
$\displaystyle (b,c) \in G^{(2)}$ then $\displaystyle (a,bc) \in G^{(2)}$ and 
$\displaystyle (ab, c) \in G^{(2)}$ and we have $a(bc) = (ab)c$, 
\item[(ii)] $\displaystyle (a,a^{-1}) \in G^{(2)}$ and 
$\displaystyle (a^{-1},a) \in G^{(2)}$, 
\item[(iii)] if $\displaystyle (a,b) \in G^{(2)}$ then $\displaystyle a b b^{-1} = a$ and 
$\displaystyle a^{-1} a b = b$. 
\end{enumerate}
\label{defgroid2}
\end{defi}

Starting with the definition \ref{defgroid2}, we can reconstruct the objects 
from definition \ref{defgroid}. The set $X= Ob(G)$ is formed by all products 
$\displaystyle a^{-1} a$, $a \in G$. For any $a \in G$ we let $\alpha(a) = 
a^{-1} a$ and $\omega(a) = a a^{-1}$. The identity section is just the identity 
function on $X$. 

\subsection{Notations.} A groupoid is denoted either by 
$\mathcal{G} = (X,G, \omega, \alpha, e, m, \, inv)$, or by 
$(G, m , inv)$. In the second case we shall use the notation 
$\displaystyle X = Ob(G) = \left\{ a^{-1} a \mbox{ : } a \in G \right\}$. 
In most of this paper  
we shall simply denote a groupoid $(G, m , inv)$ by $G$.

\begin{defi}
The transformation  $(f,F): G \rightarrow G'$ is a {\bf morphism of groupoids} defined from 
$\mathcal{G} = 
(X,G, \omega, \alpha, e, m, \, inv)$ to  $\mathcal{G}' = 
(X',G', \omega', \alpha', e', m', \, inv')$ is a pair of maps: $f: X \rightarrow 
X'$ and $F: G \rightarrow G'$ which commutes with the structure, that is: 
$f \, \omega \, = \, \omega' \, F$, $f \, \alpha \, = \, \alpha' \, F$, 
$F \, e \, = \, e' \, f$, $F \, inv \, = \, inv' \, F$ and $F$ is a morphism of 
operations, from the operation $m$ to operation $m'$.
\end{defi}

\begin{defi}
A {\bf Hausdorff topological groupoid} is a groupoid $G$ which is also a Hausdorff 
topological space, such that inversion is continuous and the multiplication is
continuous with respect to the topology on $\displaystyle G^{(2)}$ induced by the
product topology on $\displaystyle G^{2}$. 
\end{defi}

We denote by $\displaystyle dif: G \times_{\alpha} G \rightarrow G$ the 
{\bf difference function}: 
$$dif(g,h) = g h^{-1} \quad \forall (g,h) \in G \times G \quad \alpha(g) =
\alpha(h)$$

\subsection{Norms}

We shall consider the
convergence of nets $\displaystyle (a_{\varepsilon})$ of arrows, with 
$\varepsilon \in I$ a parameter in a directed set $I$. In this paper the 
most encountered directed set $I$ will be $(0,+\infty)$. 

\begin{defi}
A {\bf normed groupoid} $(G, d)$ is a groupoid 
$\mathcal{G} = (X,G, \omega, \alpha, e, m, \, inv)$ with a {\bf norm} function 
$d: G \rightarrow [0, + \infty)$, such that: 
\begin{enumerate}
\item[(i)] $d(g) = 0$ if and only if there is a $x \in X$ with $g = e(x)$, 
\item[(ii)] for any $\displaystyle (g,h) \in G^{(2)}$ we have 
$d(gh) \, \leq \, d(g) + d(h)$,  
\item[(iii)] for any $g \in G$ we have $d(inv(g)) \, = \, d(g)$. 
\end{enumerate}
A {\bf norm} $d$ {\bf is separable} if it satisfies the property: 
\begin{enumerate}
\item[(iv)] if there is a  net $\displaystyle (a_{\varepsilon}) \subset G$ such that 
for any $n \in \mathbb{N}$ $\displaystyle \alpha(a_{\varepsilon}) = x$, $\omega(a_{\varepsilon}) =
y$ and  $\displaystyle \lim_{\varepsilon \in I} d(a_{\varepsilon}) = 0$ then $x = y$.  
\end{enumerate}
\label{defnorm}
\end{defi}

\subsection{Other groupoids associated to a normed groupoid}
Let $(G,m,inv,d)$ be a normed groupoid and $dif$ its difference function. The
norm $d$ composed with the function $dif$ gives a new function $\displaystyle \tilde{d}$:  
$$\tilde{d}: G \times_{\alpha} G \rightarrow [0,+\infty) \quad , \quad \tilde{d}(g,h) =
d \, dif(g,h) = d(gh^{-1})$$ 
which induces a distance on 
the space $\displaystyle \alpha^{-1}(x)$, for any $x \in X = Ob(G)$: 
$$d_{x}: \alpha^{-1}(x) \times \alpha^{-1}(x) \rightarrow [0,+\infty) \quad ,
\quad d_{x}(g,h) = \tilde{d}d(g,h) = d(gh^{-1})$$

\begin{defi}
The {\bf metric groupoid} $\displaystyle G_{m}$  associated to the normed
groupoid $G$ is the following metric groupoid: 
\begin{enumerate}
\item[-]  the objects of $\displaystyle G_{m}$  are  the metric spaces 
$\displaystyle (\alpha^{-1}(x), d_{x})$, with $x \in Ob(G)$; 
\item[-] the arrows are  right translations 
$$\displaystyle R_{u}: (\alpha^{-1}(\omega(u)), d_{\omega(u)}) 
\rightarrow (\alpha^{-1}(\alpha(u)), d_{\alpha(u)}) \quad , \quad R_{u}(g) = gu 
$$
\item[-] the multiplication of arrows is the composition of functions; 
\item[-] the norm is defined by:  $\displaystyle 
d_{m}(R_{u}) = d_{\alpha(u)}(\alpha(u), u) = d(u)$. 
\end{enumerate}
\label{dmgru}
\end{defi} 
Remark that arrows in the metric groupoid $\displaystyle G_{m}$ are isometries. 
It is also clear that $\displaystyle G_{m}$ is isomorphic with $G$ by the morphism 
$\displaystyle u \in G \mapsto R_{u}^{-1}$.  

\begin{defi}
The {\bf $\alpha$-double groupoid} $\displaystyle G\times_{\alpha}G$ associated
to $G$ is another way to assembly the metric spaces $\displaystyle
\alpha^{-1}(x)$, $x \in Ob(G)$, into a groupoid. The definition of this groupoid
is: 
\begin{enumerate}
\item[-] the arrows are $\displaystyle G \times_{\alpha}G = \bigcup_{x \in
Ob(G)} \alpha^{-1}(x) \times \alpha^{-1}(x)$; 
\item[-] the composition of arrows is: $(g,h) (h,l) = (g,l)$, 
the inverse is $\displaystyle (g,h)^{-1} = (h,g)$,  therefore 
as a groupoid $\displaystyle G\times_{\alpha}G$ is just the union of trivial
groupoids over $\displaystyle
\alpha^{-1}(x)$, $x \in Ob(G)$; 
\item[-] it follows that $\displaystyle Ob(G\times_{\alpha}G) = 
\left\{ (g,g) \mbox{ : } g \in G \right\}$ and the induced $\alpha$ and 
$\omega$ maps are : $\displaystyle \tilde{\alpha}(g,h) = (h,h)$ and 
$ \displaystyle \tilde{\omega}(g,h) = (g,g)$, for any $g,h \in G$ with 
$\alpha(g) = \alpha(h)$; 
\item[-] the norm is the function $\displaystyle \tilde{d}$. 
\end{enumerate}
\label{ddoub}
\end{defi}
This groupoid has the property that $dif$ is a morphism of normed groupoids.

Finally, suppose that for any $x ,y \in Ob(G)$ there is $g \in G$ such that $\alpha(g) =
x$ and $\omega(g) = y$. Then any separable norm $d$ on $G$ 
induces a distance  on $X = Ob(G)$, by the 
formula: 
$$d_{ob}(x,y) = \inf \left\{ d(g) \mbox{ : } \alpha(g) = x , \omega(g) = y \right\}$$
If the groupoid is not connected by arrows then $\displaystyle d_{ob}$ may take 
the value $+ \infty$ and the space $X$ decomposes into a disjoint union of
metric spaces.

\subsection{Notions of convergence}
Any norm $d$ on a groupoid $G$ induces three notions  of
convergence on the set of arrows $G$. 


\begin{defi}
 A net of arrows  $\displaystyle (a_{\varepsilon}) $ 
 {\bf simply converges} to the arrow $a \in G$ (we write $\displaystyle a_{\varepsilon}
\rightarrow a$) if: 
\begin{enumerate}
\item[(i)] for any $\displaystyle \varepsilon \in I$ there are elements 
$\displaystyle g_{\varepsilon}, h_{\varepsilon} \in G$ such that $\displaystyle h_{\varepsilon} a_{\varepsilon} g_{\varepsilon} = 
a$, 
\item[(ii)] we have $\displaystyle \lim_{\varepsilon \in I} d(g_{\varepsilon}) =
0$ and $\displaystyle \lim_{\varepsilon \in I} d(h_{\varepsilon}) = 0$. 
\end{enumerate}

 A net of arrows  $\displaystyle (a_{\varepsilon})$ {\bf left-converges} 
 to the arrow $a \in G$ (we write 
$\displaystyle a_{\varepsilon} \lrightarrow a$) if for all  $i \in I$ we have 
$\displaystyle (a_{\varepsilon}^{-1}, a) \in G^{(2)}$ and 
moreover $\displaystyle \lim_{\varepsilon \in I} d(a_{\varepsilon}^{-1} a) = 0$. 

 A net of arrows  $\displaystyle (a_{\varepsilon})$ {\bf right-converges} 
 to the arrow $a \in G$ (we write 
$\displaystyle a_{\varepsilon} \rrightarrow a$) if for all  $i \in I$ we have 
 $\displaystyle (a_{\varepsilon}, a^{-1}) \in G^{(2)}$ and 
moreover $\displaystyle \lim_{\varepsilon \in I} d(a_{\varepsilon} a^{-1}) = 0$. 
\label{defconv}
\end{defi}

It is clear that if $\displaystyle a_{\varepsilon} \lrightarrow a$ or 
$\displaystyle a_{\varepsilon} \rrightarrow a$ then 
$\displaystyle a_{\varepsilon} \rightarrow a$. 

Right-convergence of $\displaystyle a_{\varepsilon}$ to $a$ is just convergence 
of $\displaystyle a_{\varepsilon}$ to $a$ in the distance $\displaystyle
d_{\alpha(a)}$, that is $\displaystyle \lim_{\varepsilon \in I}
d_{\alpha(a)}(a_{\varepsilon}, a)= 0$. 

Left-convergence of $\displaystyle a_{\varepsilon}$ to $a$ is just convergence 
of $\displaystyle a_{\varepsilon}^{-1}$ to $a^{-1}$ in the distance 
$\displaystyle
d_{\omega(a)}$, that is $\displaystyle \lim_{\varepsilon \in I}
d_{\omega(a)}(a_{\varepsilon}^{-1}, a^{-1})= 0$.

\begin{prop}
\label{propconv}
Let $G$ be a groupoid with a norm $d$.  
\begin{enumerate}
\item[(i)] If $\displaystyle a_{\varepsilon} \lrightarrow a$ and $\displaystyle a_{\varepsilon}
\lrightarrow b$ then $a=b$. If  $\displaystyle a_{\varepsilon} \rrightarrow a$ and $\displaystyle a_{\varepsilon}
\rrightarrow b$ then $a=b$. 
\item[(ii)] The following are equivalent: 
\begin{enumerate}
\item[1.] $G$ is a Hausdorff topological groupoid with respect to the topology induced by 
the simple convergence, 
\item[2.] $d$ is a separable norm, 
\item[3.] for any net $\displaystyle (a_{\varepsilon}) $, if $\displaystyle a_{\varepsilon} 
\rightarrow a$ and $\displaystyle a_{\varepsilon}
\rightarrow b$ then $a=b$.
\item[4.] for any net $\displaystyle (a_{\varepsilon}) $, if $\displaystyle a_{\varepsilon} 
\rrightarrow a$ and $\displaystyle a_{\varepsilon}
\lrightarrow b$ then $a=b$.
\end{enumerate}
\end{enumerate}
\end{prop}

\paragraph{Proof.}
(i) We prove only the first part of the conclusion.  We can write $\displaystyle b^{-1} a = b^{-1} a_{\varepsilon} 
a_{\varepsilon}^{-1} a$, therefore 
$$d(b^{-1}a) \leq d(b^{-1} a_{\varepsilon}) + d(a_{\varepsilon}^{-1} a)$$
The right hand side of this inequality is arbitrarily small, so $\displaystyle 
d(b^{-1}a) =0$, which implies $a=b$. 

(ii) Remark that the structure maps are continuous with
respect to the topology induced by the simple convergence. We need only to prove
the uniqueness of limits. 

 3. $\Rightarrow$ 4. is trivial. In order to prove that 
4.$\Rightarrow$ 3., consider an arbitrary  net  
$\displaystyle (a_{\varepsilon}) $ such that  $\displaystyle a_{\varepsilon} 
\rightarrow a$ and $\displaystyle a_{\varepsilon}
\rightarrow b$. This means that there exist nets $\displaystyle 
(g_{\varepsilon}) , (g'_{\varepsilon}) , (h_{\varepsilon}) , 
(h'_{\varepsilon}) $ such that 
$\displaystyle h_{\varepsilon} a_{\varepsilon} g_{\varepsilon} = a$, 
$\displaystyle h'_{\varepsilon} a_{\varepsilon} g'_{\varepsilon} = b$ 
and $\displaystyle \lim_{i \in I} \left( d(g_{\varepsilon}) + 
d(g'_{\varepsilon}) + d(h_{\varepsilon}) + d(h'_{\varepsilon}) \right) = 0$. 
Let $\displaystyle g"_{\varepsilon} = g_{\varepsilon}^{-1} g'_{\varepsilon}$ 
and $\displaystyle h"_{\varepsilon} = h'_{\varepsilon} h_{\varepsilon}^{-1}$. 
We have then $\displaystyle b = h"_{\varepsilon} a g"_{\varepsilon}$ and 
$\displaystyle  \lim_{i \in I} \left( d(g"_{\varepsilon}) + d(h"_{\varepsilon})
\right) = 0$. Then $\displaystyle h"_{\varepsilon} a \lrightarrow b$ and 
$\displaystyle h"_{\varepsilon}a \rrightarrow a$. We deduce that $a=b$. 

1.$\Leftrightarrow$ 3. is trivial. So is 3. $\Rightarrow$ 2. 
We finish the proof by showing that 2. $\Rightarrow$ 3. 
By a reasoning made previously, it is enough to prove that: 
if $ b = h_{\varepsilon} a g_{\varepsilon}$ and 
$\displaystyle \lim_{i \in I} \left( d(g_{\varepsilon}) + d(h_{\varepsilon}) 
\right) = 0$ then $a=b$. Because $d$ is separable it follows that 
$\alpha(a) = \alpha(b)$ and $\omega(a) = \omega(b)$. We have then 
$\displaystyle a^{-1} b = a^{-1} h_{\varepsilon} a g_{\varepsilon}$, therefore 
$$d(a^{-1} b) \leq d(a^{-1} h_{\varepsilon} a) + d(g_{\varepsilon})$$
 The  norm $d$ induces a left invariant  distance on the vertex group of all 
 arrows $g$ such that $\alpha(g) = \omega(g) = \alpha(a)$. This distance is 
 obviously continuous with respect to the simple convergence  in the group. 
 The net $\displaystyle a^{-1} h_{\varepsilon} a$ simply converges to 
 $\alpha(a)$ by the continuity of the multiplication (indeed, 
 $\displaystyle h_{\varepsilon}$ simply converges to $\alpha(a)$). 
 Therefore  $\displaystyle \lim_{i \in I} d(a^{-1} h_{\varepsilon} a) = 0$.
It follows that $\displaystyle d(a^{-1}b)$ is arbitrarily small, therefore $a=b$. 
\hfill $\quad \square$

\subsection{Families of seminorms}
Instead of a norm we may use  families of seminorms.  

\begin{defi}
A {\bf family of seminorms}  on a groupoid $G$ 
is a family $S$ of functions $\rho: G \rightarrow [0,+\infty)$ with the properties: 
\begin{enumerate}
\item[(i)] for any $x \in X$ and $\rho \in S$ we have $\rho(e(x)) = 0$; 
 if $\rho(g) = 0$ for any $\rho \in S$ then   there is $x \in X$ such that 
$g = e(x)$, 
\item[(ii)] for any $\rho \in S$ and $(g,h) \in G^{(2)}$ we have 
$\displaystyle \rho(gh) \leq \rho(g) + \rho(h)$, 
\item[(iii)] for any $\rho \in S$ and $g \in G$ we have $\rho(inv(g)) = \rho(g)$. 
\end{enumerate}
A groupoid $G$ endowed with a  family of seminorms $S$ is called a {\bf 
  seminormed groupoid}. 

A   family of seminorms $S$ is separable if it satisfies the property: 
\begin{enumerate}
\item[(iv)] if there is a  net $\displaystyle (a_{\varepsilon})  \subset G$ such that 
for any $n \in \mathbb{N}$ $\displaystyle \alpha(a_{\varepsilon}) = x$, $\omega(a_{\varepsilon}) =
y$ and  for any $\rho \in S$ we have 
$\displaystyle \lim_{i \in I} \rho(a_{\varepsilon}) = 0$ then $x = y$. 
\end{enumerate}
\label{defseminormedg}
\end{defi}

Families of morphisms  induce   families of seminorms. 

\begin{defi}
Let $G$ be a groupoid and $(H,d)$ be a normed groupoid. A $(H,d)$ {\bf family of
morphisms} is a set  $L$ of
morphisms from $G$ to $H$ such that  for any 
$g \in G$ there is $A \in L$ with $A(g) \not \in Ob(H)$. 
\end{defi}

The following proposition has a straightforward proof which we omit. 

\begin{prop}
Let  $G$ be a groupoid, $(H,d)$ be a normed groupoid and $L$ a    $(H,d)$
family of morphisms. Then the set 
$$d(L) = \left\{ d \, A  \mbox{ : } A \in L \right\}$$
is a family of seminorms. 
\end{prop}

Definition \ref{defconv} can be modified for the case of   families of
seminorms. 

\begin{defi}
Let $(G,S)$ be a   semi-normed groupoid.  A net of arrows  
$\displaystyle (a_{\varepsilon})$ {\bf simply converges} to the arrow 
$a \in G$ (we write $\displaystyle a_{\varepsilon} \rightarrow a$) if: 
\begin{enumerate}
\item[(i)] for any $\displaystyle i \in I$ there are elements 
$\displaystyle g_{\varepsilon}, h_{\varepsilon} \in G$ such that 
$\displaystyle h_{\varepsilon} a_{\varepsilon} g_{\varepsilon} = a$, 
\item[(ii)] for any $\rho \in S$  we have $\displaystyle 
\lim_{i \in I} \left( \rho(g_{\varepsilon})+ \rho(h_{\varepsilon}) \right) = 0$. 
\end{enumerate}

 A net of arrows  $\displaystyle (a_{\varepsilon})$ {\bf left-converges} 
 to the arrow $a \in G$ (we write 
 $\displaystyle a_{\varepsilon} \lrightarrow a$) if for all 
 $i \in I$ we have $\displaystyle (a_{\varepsilon}^{-1}, a) \in G^{(2)}$ 
 and moreover for any $\rho \in S$ we have 
 $\displaystyle \lim_{i \in I} \rho(a_{\varepsilon}^{-1} a) = 0$. 

 A net of arrows  $\displaystyle (a_{\varepsilon})$ {\bf right-converges} 
 to the arrow $a \in G$ (we write  
 $\displaystyle a_{\varepsilon} \rrightarrow a$) if for all $i \in I$ we have 
 $\displaystyle (a_{\varepsilon}, a^{-1}) \in G^{(2)}$ and 
moreover  for any $\rho \in S$ we have $\displaystyle \lim_{i \in I} 
\rho(a_{\varepsilon} a^{-1}) = 0$. 
\label{defconv2}
\end{defi}

With these slight modifications,  the proposition \ref{propconv} still holds
true. This is visible from the examination of its proof. 

Let us finally remark that if $(G, dL)$ is a   seminormed groupoid, where 
$L$ is a   $(H,d)$ family of morphisms, then a net  $\displaystyle 
(a_{\varepsilon}) \in G$ converges (simply, left or right) to $a \in G$ if and 
only if for any $A \in L$ the net $\displaystyle (A(a_{\varepsilon}))$ 
respectively converges in $(H,d)$. 

\subsection{Uniform convergence on bounded sets}
\label{secdefconv}

We shall use right-convergence, according to definition 
\ref{defconv}, but left-convergence or simple convergence could also be used. In
relation to this see for example the remark \ref{rema1}.

\begin{defi}
Let $G$ be a normed groupoid with a separable norm.  
A net  $\displaystyle (f_{\varepsilon})$  of functions 
$\displaystyle f_{\varepsilon}: G\times_{\alpha}G\rightarrow G$ {\bf uniformly 
converges on bounded sets} to the function $\displaystyle f: G\times_{\alpha}G\rightarrow G$ (in the sense
of the left convergence) if: 
\begin{enumerate}
\item[(i)] for any $\varepsilon > 0$ and  $\displaystyle (h,g) \in G^{(2)}$  we
have $\displaystyle \alpha(f_{\varepsilon}(h,g)) = \alpha(f(h,g))$, 
\item[(ii)] for any $\lambda, \mu > 0$ there is $\varepsilon(\lambda,\mu) > 0$ 
such that for any $\varepsilon \in \Gamma$, $\mid\varepsilon\mid \leq \varepsilon(\lambda,\mu)$ and any 
$\displaystyle (h,g) \in G^{(2)}$ with $d(h) \leq \lambda$, $d(g) \leq \lambda$,
we have: 
$$d\left( f_{\varepsilon}(h,g) inv(f(h,g)) \right) \leq \mu$$
\end{enumerate}
In the case of a groupoid $G$ with a separable  family of seminorms $S$, the definition of
uniform convergence is the same, excepting the modification of (ii) above into: 
for  any $\lambda, \mu > 0$ and any seminorm $\rho \in S$ there is 
$\varepsilon(\lambda,\mu, \rho) > 0$ such that for any $\varepsilon \in \Gamma$, 
$\mid\varepsilon\mid \leq \varepsilon(\lambda,\mu)$ and any $\displaystyle (h,g) \in G^{(2)}$ 
with $\rho(h) \leq \lambda$, $\rho(g) \leq \lambda$,
we have: 
$$\rho\left(  f_{\varepsilon}(h,g) inv(f(h,g))\right) \leq \mu$$

Similarly, in a normed groupoid with a separable norm $d$, 
the uniform convergence on bounded sets of a net of functions $\displaystyle f_{\varepsilon}: G
\rightarrow \mathbb{R}$ to $f:G \rightarrow \mathbb{R}$ means that 
for any $\lambda, \mu > 0$ there is $\varepsilon(\lambda,\mu) > 0$ 
such that for any $\varepsilon \in \Gamma$, $\mid\varepsilon\mid \leq \varepsilon(\lambda,\mu)$ and any 
$\displaystyle g \in G$ with $d(g) \leq \lambda$ 
we have: $\mid f_{\varepsilon}(g) - f(g)\mid \leq \mu$. 
\label{defcon}
\end{defi}

\section{Normed categories}

Sometimes it is interesting to work with normed categories, instead of normed
groupoids. Briefly said, a normed category is a small category endowed with an
involutionary antimorphism (an "inverse") and with a norm function 
defined on arrows.

\begin{defi}
A {\bf normed category} $(G, d)$ is a small category $\mathcal{C}$ 
with inverses:    the set of  
objects is $X = Obj(\mathcal{C})$, 
 the set of arrows is identified with $\mathcal{C}$; there are 
two functions $s, t: \mathcal{C} \rightarrow X$, named {\bf source} and {\bf
target}, a {\bf multiplication}  
$m : \mathcal{C}^{(2)} \rightarrow \mathcal{C}$ (notation: $m(g,h) = gh$ or 
$m(g,h) = g \circ h$), 
where 
$$\mathcal{C}^{(2)} \, = \, \left\{ (g,h) \in \mathcal{C}^{2} \mbox{ : } s(g) =
t(h) \right\}$$
 a {\bf inversion map} $inv: \mathcal{C} \rightarrow \mathcal{C}$ (notation 
 $\displaystyle inv(g) = g^{-1}$) and a {\bf
 norm} $d : \mathcal{C} \rightarrow [0,+\infty)$, which
satisfy the following conditions: 
\begin{enumerate}
\item[(i)] $s(gh) = s(h)$, $t(gh) = t(g)$, for any $\displaystyle 
(g,h) \in \mathcal{C}^{(2)}$, 
\item[(ii)] the multiplication is associative $(gh)k = g(hk)$ 
\item[(iii)] the inverse is an  involution and a antimorphism: 
  $\displaystyle s(g^{-1}) = t(g)$, $\displaystyle \left( gh \right)^{-1} = h^{-1}
  g^{-1}$. 
\item[(iv)] $d(g) = 0$ if and only if there is $h \in \mathcal{C}$ such that 
$\displaystyle g = h^{-1} h$,
\item[(v)] for any $\displaystyle (g,h) \in \mathcal{C}^{(2)}$ we have 
$d(gh) \, \leq \, d(g) + d(h)$,  
\item[(vi)] for any $g \in \mathcal{C}$ we have $d(inv(g)) \, = \, d(g)$. 
\end{enumerate}
Let $\displaystyle \alpha, \omega : \mathcal{C} \rightarrow \mathcal{C}$ be
defined by $\displaystyle \alpha(g) = g^{-1} g$, $\displaystyle \omega(g) = g
g^{-1}$. 
The {\bf norm} $d$ {\bf is separable} if it satisfies the property: 
\begin{enumerate}
\item[(vii)] if there is a  net $\displaystyle (a_{\varepsilon}) \subset
\mathcal{C}$ such that 
for any $n \in \mathbb{N}$ $\displaystyle \alpha(a_{\varepsilon}) = x$, 
$\omega(a_{\varepsilon}) =
y$ and  $\displaystyle \lim_{\varepsilon \in I} d(a_{\varepsilon}) = 0$ then 
$x = y$.  
\end{enumerate}
\label{defcatnorm}
\end{defi}

Seminormed categories are defined further, by making a slight modification 
of definition \ref{defseminormedg}. 

\begin{defi}
Let $\mathcal{C}$ be a category with inverses (which satisfies (i)--(iv) 
definition \ref{defcatnorm}). Let $\displaystyle X = \alpha(\mathcal{C}$. 

A {\bf family of seminorms}  on a category with inverses $\mathcal{C}$ 
is a family $S$ of functions $\rho: \mathcal{C} \rightarrow [0,+\infty)$ 
with the properties: 
\begin{enumerate}
\item[(iv)'] for any $x \in X$ and $\rho \in S$ we have $\rho(x) = 0$; 
 if $\rho(g) = 0$ for any $\rho \in S$ then    $g \in X$ , 
\item[(v)'] for any $\rho \in S$ and $(g,h) \in \mathcal{C}^{(2)}$ we have 
$\displaystyle \rho(gh) \leq \rho(g) + \rho(h)$, 
\item[(vi)'] for any $\rho \in S$ and $g \in \mathcal{C}$ we have $\rho(inv(g)) = \rho(g)$. 
\end{enumerate}
A category $\mathcal{C}$ with inverses endowed with a  family of seminorms 
$S$ is called a {\bf   seminormed category}. 

A   family of seminorms $S$ is separable if it satisfies the property: 
\begin{enumerate}
\item[(vii)'] if there is a  net $\displaystyle (a_{\varepsilon})  \subset
\mathcal{C}$ 
such that  for any $i \in I$ 
$\displaystyle \alpha(a_{\varepsilon}) = x$, $\omega(a_{\varepsilon}) =
y$ and  for any $\rho \in S$ we have 
$\displaystyle \lim_{i \in I} \rho(a_{\varepsilon}) = 0$ then $x = y$. 
\end{enumerate}
\label{defseminormedc}
\end{defi}

All considerations made before concerning convergence for normed or seminormed
separable groupoids, extend without effort to normed categories, or to
seminormed separable categories. A little bit of care is needed though:
everywhere we should replace source and target maps (of groupoids) by 
the (algebraically defined) $\displaystyle \alpha(g) = g^{-1}g$ and $\omega(g) =
g g^{-1}$. 

\section{Examples of  normed groupoids and normed categories}

We give several examples of normed groupoids which will be of interest later in
this paper. 

\subsection{Metric spaces}
\label{smetricspace}

Let $(X,d)$ be a metric space. We form the {\bf normed trivial groupoid} $(G,d)$ over
$X$: 
\begin{enumerate}
\item[-] the set of arrows is $G = X \times X$ and the multiplication is 
$$(x,y) (y,z) = (x,z)$$
Therefore we have $\alpha(x,y) = y$, $\omega(x,y) = x$, $e(x) = (x,x)$, 
$\displaystyle (x,y)^{-1} = (y,x)$. 
\item[-] the norm is just the distance function $d: G \rightarrow [0,+\infty)$.
\end{enumerate}

It is easy to see that if $(X \times X, d)$ is a normed trivial groupoid over $X$
then $(X,d)$ is a metric space. 

\subsubsection{Associated groupoids.} The metric groupoid $\displaystyle (X\times
X)_{m}$ can be described as the groupoid with objects pointed metric spaces 
 $(X,d, x)$, $x \in X$, arrows $\displaystyle R_{(x,y)} : (X,d, x)
 \rightarrow (X,d,y)$, $\displaystyle  R_{(x,y)}(z) = z$, and norm 
 $\displaystyle d_{m}\left( R_{(x,y)} \right) = d(x,y)$. The $\alpha$-double 
 groupoid $(X\times X) \times_{\alpha} (X \times X)$ can be described as the 
 groupoid with arrows $X \times X \times X$ , composition 
 $(x,y,z)(y,v,z) = (x,v,z)$, inverse $\displaystyle (x, y, z)^{-1} = (y,x,z)$ 
 and norm $\displaystyle \tilde{d}(x,y,z) = d(x,y)$. 
   
\subsubsection{Convergence.} Remark first that $d$ is a separable norm, according to definition 
\ref{defnorm} (iv). Indeed, for
any $x,y \in X$ there is only one arrow $a \in X \times X$ such that $\alpha(a)
= x$, $\omega(a) = y$, namely the arrow $a= (y,x)$. Any net 
$\displaystyle (a_{\varepsilon}) $ with $\displaystyle \alpha(a_{\varepsilon})
= x$, $\omega(a_{\varepsilon}) = y$ is the constant net $\displaystyle a_{\varepsilon} =
(y,x)$. If $\displaystyle \lim_{n \rightarrow \infty} d(a_{\varepsilon}) = 0$ then 
$d(y,x) = 0$, therefore $x=y$. We deduce from proposition \ref{propconv} that 
we have only one interesting notion of convergence, which is simple convergence.

In the particular case of normed trivial groupoids the definition \ref{defconv} 
of simple convergence becomes: a net $\displaystyle (x_{\varepsilon}, y_{\varepsilon})  \subset (X
\times X$ simply converges to $(x,y)$ if we have 
$$\displaystyle \lim_{n \rightarrow
\infty} \left( d(x,x_{\varepsilon}) + d(y_{\varepsilon}, y) \right) = 0$$ 
 that is if the nets $\displaystyle 
x_{\varepsilon}, y_{\varepsilon}$  converge respectively to $x, y$. Indeed this is coming from the
fact that for any $n \in \mathbb{N}$ there are unique $\displaystyle h_{\varepsilon},
g_{\varepsilon} \in X \times X$ such that $\displaystyle h_{\varepsilon} (x_{\varepsilon}, y_{\varepsilon}) g_{\varepsilon} = 
(x,y)$. These are $\displaystyle h_{\varepsilon} = (x,x_{\varepsilon})$ and $\displaystyle 
g_{\varepsilon} = (y_{\varepsilon}, y)$.

\subsubsection{Nice families of seminorms on metric spaces.} Let $X$ be a non empty 
set, let 
$(Y,d)$ be a metric spaces and  $\displaystyle 
(Y^{2}, d)$ its associated normed trivial groupoid. Any function 
$f: X \rightarrow Y$ induces a morphism $\displaystyle \bar{f}$ from the trivial
groupoid $\displaystyle X^{2}$ to $\displaystyle Y^{2}$ by 
$\displaystyle \bar{f}(x,y) = (f(x), f(y))$. Any  family $ \mathcal{L}$ 
of  functions from 
$X$ to $Y$  with the separation property: for any $x, y \in X$ $x \not = y$
there is $f \in \mathcal{L}$ with $f(x) \not = f(y)$, 
gives us a   $\displaystyle (Y^{2}, d)$ family of morphisms, which in turn induces a   family of
seminorms on $\displaystyle X^{2}$.

\subsection{Normed groupoids from $\alpha$-double groupoids}

We can construct normed groupoids starting from definition \ref{ddoub} 
of $\alpha$-double groupoids. 

\begin{prop}
Let $(G,d)$ be a groupoid and $\displaystyle (G\times_{\alpha}G,\tilde{d})$ the 
associated $\alpha$-double groupoid. Then for any $\displaystyle (g,h) \in 
G\times_{\alpha}G$ and for any $u \in G$ with $\omega(u) = \alpha(g) =
\alpha(h)$ we have 
\begin{equation}
d_{\omega(u)}(g,h) = d_{\alpha(u)}(gu, hu)
\label{leftdis}
\end{equation}
Conversely, suppose that $G$ is a groupoid and that for any $x \in Ob(G)$ we
have a distance $\displaystyle d_{x}: \alpha^{-1}(x) \times \alpha^{-1}(x)
\rightarrow [0,+\infty)$. If (\ref{leftdis}) is true for any 
$\displaystyle (g,h) \in G\times_{\alpha}G$ and for any $u \in G$ with 
$\omega(u) = \alpha(g) = \alpha(h)$ then 
$$d(g) = d_{\alpha(g)}(g,\alpha(g)) \quad \mbox{ and } \quad \tilde{d}(g,h) = 
d_{\alpha(g)}(g,h)$$
define a norm on $G$ such that $\displaystyle (G\times_{\alpha}G, \tilde{d})$ 
is the associated $\alpha$-double groupoid. 
\label{pdoub}
\end{prop}

\begin{rk}
Therefore any normed groupoid $(G,d)$ can be seen as the  bundle  of 
metric spaces $\displaystyle \alpha: G \rightarrow Ob(G)$, such that 
(a) each fiber $\displaystyle \alpha^{-1}(x)$ has a distance $\displaystyle 
d_{x}$, and (b) the distances $\displaystyle d_{x}$ are right invariant 
with respect to the groupoid composition, in the sense of relation
(\ref{leftdis}). 
\end{rk}

\paragraph{Proof.} 
For the first implication remark that $\displaystyle (gu, hu) \in
G\times_{\alpha}G$. Moreover let $g' =gu$, $h' = hu$. Then $\displaystyle 
g' \left(h'\right)^{-1} = g h^{-1}$, therefore 
$$d_{\omega(u)}(g,h) = d(gh^{-1}) = d_{\alpha(u)}(gu, hu)$$
For the converse implication, we have to prove that if $\displaystyle 
g' \left(h'\right)^{-1} = g h^{-1}$ then $\displaystyle d(gh^{-1}) = 
d(g' \left(h'\right)^{-1})$, with $d$ defined as in the formulation of the 
proposition. This is easy: Let $u = \left(h'\right)^{-1}) h$, then $g = g'u$, 
$h = h'u$ and (\ref{pdoub}) implies the desired equality. The verification that
$d$ is indeed a norm on $G$ is straightforward, as well as the fact that 
$\displaystyle \tilde{d}$ is the induced norm on $\displaystyle
G\times_{\alpha}G$. \quad $\square$

\subsection{Group actions}
\label{saction}

Let $G$ be a group with neutral element $e$, 
 which acts from the left on the space $X$. Associated with
this is the {\bf action groupoid} $G \ltimes X$ over $X$. 
The action groupoid is defined as: the set of arrows is $X \times G$ and the multiplication is 
$$(g(x),h) (x,g) = (x, hg)$$
Therefore $\displaystyle \alpha(x,g) = x$, $\omega(x,g) = g(x)$, 
$e(x) = (x,e)$, $\displaystyle (x,g)^{-1} = (g(x), g^{-1})$,

As a particular case of definition \ref{defnorm}, a 
 {\bf normed action groupoid} is an action groupoid $G \ltimes X$ endowed with 
a norm function $d: X \times G \rightarrow [0,+ \infty)$ with the properties: 
\begin{enumerate}
\item[(i)] $d(x,g) = 0$ if and only if $g=e$, 
\item[(ii)] $\displaystyle d(g(x), g^{-1}) = d(x,g)$, 
\item[(iii)] $d(x, hg) \leq d(x,g) + d(g(x), h)$. 
\end{enumerate}

Remark that the norm function is no longer a distance function. 
In the case of a free action (if $g(x) = x$ for some $x \in X$ then $g = e$) we
may obtain a norm function from a distance function on $X$. Indeed, let 
$d': X \times X \rightarrow [0,+\infty)$ be a distance. Define then 
$\displaystyle d: X \times G \rightarrow [0,+\infty)$ by 
$$\bar{d}(x,g) = d'(g(x), x)$$
Then $\displaystyle (G \ltimes X, \bar{d})$ is a normed action groupoid. 

\subsubsection{Associated groupoids.} The associated $\alpha$-double groupoid can be seen as 
$X \times G \times G$, with composition $(x,g,h)(x,h,l) = (x,g,l)$ and 
inverse $\displaystyle (x,g,h)^{-1} = (x,h,g)$. For any $x \in X$ we have 
a distance $\displaystyle d_{x}: G \times G \rightarrow [0,+ \infty)$, defined
by 
$$d_{x}(g,h) = d(h(x), gh^{-1})$$
Conversely, according to proposition \ref{pdoub} and relation (\ref{leftdis}), 
a norm on a action groupoid can be constructed from a function $\displaystyle x \in X \mapsto d_{x}$ which 
associates to any $x \in X$ a distance $\displaystyle d_{x}$ on G, such that 
for any $x \in X$ and $u,g,h \in G$ we have 
$$d_{u(x)}(g,h) = d_{x}(gu, hu)$$
In this case we can define the norm on the action groupoid by 
 $\displaystyle d(x,g) = d_{x}(g,e)$. 

A particular case is $X = \left\{x \right\}$, when a normed action groupoid is
just a group endowed with a right invariant distance.  

\subsubsection{Convergence.} The norm $d$ is separable if the following condition is
satisfied: for any $x, y \in X$ and any net 
$\displaystyle g_{\varepsilon} \in G$ with the property $\displaystyle
g_{\varepsilon}(x) = y$ for all $\varepsilon$, if $\displaystyle
\lim_{\varepsilon \in I} d_{x}(g_{\varepsilon}, e) = 0$ then $x = y$.

\subsection{Groupoids actions}

Let $\displaystyle G$  and $M$ be two groupoids. We denote 
by $Aut(M)$ the groupoid which has as objects sub-groupoids of $M$ and 
invertible morphisms between sub-groupoids of $M$ as arrows. A groupoid action 
of $G$ on $M$ is just a morphism $F$ of groupoids from $G$ to $Aut(M)$. In fewer
words, for any $g \in G$ let $F(g)$ be the associated morphism of 
sub-groupoids, defined from the sub-groupoid denoted by $dom \, g$ to the
sub-groupoid denoted by $im \, g$. For any $x \in dom \, g$ we use the notation 
 $g . x = F(g)(x)$. Compositions in $G$ and in $M$ are denoted 
 multiplicatively. Let  $G \ltimes M$ be the set 
 $$ G \ltimes M \, = \, \left\{ (x, g) \mbox{ : } x \in dom \, g \right\}$$ 
 The action of $G$ on $M$ satisfies the following conditions: 
 \begin{enumerate}
 \item[-] for any $u, v \in G$ such that $\alpha(v) = \omega(u)$ we have 
 $dom \, u \, = \, dom \, vu$, $im \, u \, = \, dom \, v$ and for any 
 $x \in dom \, u$ we have $v . ( u . x) \, = \, (vu) . x$; 
 \item[-] for any $u \in G$ and $x, y \in dom \, u$ we have $u . (xy) \, = \, 
 (u . x) (u . y)$. 
 \end{enumerate}
 
 Any groupoid action induces a groupoid structure on $G \ltimes M$, by the
 composition law $(g . x , h) (x , g) \, = \, (x , hg)$. 
 
 At a closer look we may notice an example of a groupoid action in 
 proposition \ref{pdoub}.  
 Indeed, let $G$ be a groupoid and $\displaystyle G \times_{\alpha} G$ 
 the associated $\alpha$-double groupoid. Then $G$ acts on 
 $\displaystyle G \times_{\alpha} G$ by $u . (g,h) \, = \, (g u^{-1} , h
 u^{-1})$, for any $u \in G$ and any $\displaystyle (g,h) \in G \times_{\alpha}
 G$ such that $\alpha(u) = \alpha(g) = \alpha(h)$. Therefore 
 $\displaystyle dom \, u \, = \, \left( \alpha^{-1}(\alpha(u)) \right)^{2}$ and 
 the associated action groupoid is 
 $$G \ltimes (G \times_{\alpha} G) \, = \, \left\{ (g,h, u) \mbox{ : } 
 \alpha(g) = \alpha(h) = \alpha(u) \right\}$$ 
 with multiplication defined by 
 $$(g u^{-1} , h u^{-1} , v) (g , h, u) \, = \, (g, h , vu)$$
 Relation (\ref{leftdis}) in proposition \ref{pdoub} tells that $G$ acts on the
 normed groupoid $\displaystyle (G \times_{\alpha} G, \tilde{d})$ by 
 isometries. In general, the action groupoid induced by the action of a groupoid
 $G$ on a normed groupoid $M$ {\it by isometries} may be an object as
 interesting as a normed groupoid.

\subsection{Riemannian groupoids (not normed)}
\label{riemgrou}

The references used are Glickenstein \cite{glick} and Lott \cite{lott1} section
5. 

\begin{defi}
A smooth \'etale groupoid $G$ is Riemannian if there is a riemannian metric on
$Ob(G)$ such that the germs of the arrows (arrows are seen as local
diffeomeorphisms on $Ob(G)$) are isometries. 
\end{defi}

A Riemannian groupoid is not a normed groupoid.

\subsection{Transport category}

A good reference for this section is \cite{AGS} chapters 5 and 6. Much more details could be found, for instance, 
 in  \cite{VILLANI}.

Let $(X,d)$ be a a complete, separable metric space (Polish space). The class of Borel 
probabilities on $X$ is denoted by $\mathcal{P}(X)$. 
If $\displaystyle (X_{1},d_{1}), \, (X_{2}, d_{2})$  are Polish spaces then 
$\displaystyle (X_{1} \times X_{2}, d)$ is also a Polish space, with 
the distance $d$ defined by $\displaystyle d((x_{1}, x_{2}), (y_{1}, y_{2}) = 
d_{1}(x_{1}, y_{1}) + d_{2}(x_{2}, y_{2})$.

If $\displaystyle f: X_{1} \rightarrow X_{2}$ is a Borel map and 
$\displaystyle \mu \in \mathcal{P}(X_{1})$ then the push-forward of 
$\mu$ through $f$ is the measure $\displaystyle f \sharp \mu \in
\mathcal{P}(X_{2}$ defined by: for any $\displaystyle B\in\mathcal{B}(X_{2})$ 
$$\left(f \sharp \mu \right) (B) \, = \, \mu(f^{-1}(B))$$
In particular, the projections $\displaystyle \pi_{i}: X_{1} \times X_{2} \rightarrow X_{i}$, 
$i = 1, 2$, $\displaystyle \pi_{i}(x_{1}, x_{2}) = x_{i}$,  define push-forwards
from $\displaystyle \mathcal{P}(X_{1} \times X_{2})$ to $\mathcal{P}(X_{i})$. 

\begin{defi}
Let $X$, $Y$ be Polish spaces. A measure-valued map $\displaystyle x \in X \mapsto \mu_{x} 
\in \mathcal{P}(Y)$ is a Borel map if for any Borel set $B \in \mathcal{B}(Y)$
the function $x \mapsto \mu_{x}(B)$ is Borel. 
\end{defi}

The following is the disintegration theorem \cite{AGS} theorem 5.3.1 (see also 
references therein). (In \cite{AGS} the theorem is formulated for Radon spaces,
here we stay in the frame of Polish spaces). 

\begin{thm}
Let $X, Y$ be two Polish spaces, $\pi: X \rightarrow Y$ a Borel map and 
$\\mu \in \mathcal{P}(X)$, $\nu = \pi \sharp \mu \, \in \mathcal{P}(Y)$. Then
there exists a Borel measure-valued function $y \in Y \mapsto \mu_{y} \in \mathcal{P}(X)$, 
which is $\nu$-a.e. uniquely determined, such that: 
\begin{enumerate}
\item[(a)] for $\nu$-a.e. $y \in Y$ we have $\displaystyle 
\mu_{y}(X \setminus \pi^{-1}(x)) = 0$, 
\item[(b)] for every Borel map $f : X \rightarrow [0,+\infty]$  
$$\int_{X} f(x) \mbox{ d}\mu(x) \, = \, \int_{Y} \left( \int_{\pi^{-1}(y)} 
f(x) \mbox{ d} \mu_{y}(x) \right\} \mbox{ d} \nu(y)$$
\end{enumerate}
\label{disthm}
\end{thm}
In the particular case $\displaystyle X = X_{1} \times X_{2}$, $\displaystyle 
\pi = \pi_{i}$ and $\displaystyle Y = X_{i}$, $i = 1,2$, the disintegration
theorem implies that for any $\displaystyle \gamma \in \mathcal{P}(X_{1} \times
X_{2})$ with marginals 
$$\phi_{1}\sharp \gamma = \mu \, , \, \phi_{2} \sharp \gamma = \nu $$ 
there exist Borel measure-valued maps $\displaystyle 
x_{1} \in X_{1} \mapsto \gamma_{x_{1}} \in \mathcal{P}(X_{2})$, $\mu$-a.e. 
uniquely determined, and $\displaystyle 
x_{2} \in X_{2} \mapsto \gamma_{x_{2}} \in \mathcal{P}(X_{1})$, $\nu$-a.e. 
uniquely determined, such that for any Borel function $\displaystyle 
f: X_{1} \times X_{2} \rightarrow [0,+\infty]$ we have: 
$$\int_{X_{1} \times X_{2}} f(x_{1}, x_{2}) \mbox{ d} \gamma(x_{1}, x_{2}) \, =
\, \int_{X_{1}} \left( \int_{X_{2}} f(x_{1}, x_{2}) \mbox{
d}\gamma_{x_{1}}(x_{2}) \right) \mbox{ d}\mu(x_{1}) \, = \, $$
$$= \, 
\int_{X_{2}} \left( \int_{X_{1}} f(x_{1}, x_{2}) \mbox{
d}\gamma_{x_{2}}(x_{1}) \right) \mbox{ d}\nu(x_{2})$$

For any two probability measures $\mu$, $\nu$ on $X$ the set of all transport 
plans between $\mu$ and $\nu$ is defined as: 
$$\Pi (\mu, \nu) \, = \, \left\{ \gamma \in \mathcal{P}(X \times X) \mbox{ : } 
\phi_{1}\sharp \gamma = \mu \, , \, \phi_{2} \sharp \gamma = \nu \right\}$$

\begin{defi}
The category $Trans(X)$ of transport plans on $X$ is the category with inverses 
 which has as objects the elements of $\mathcal{P}(X)$ and the class of  arrows between 
$\mu, \nu \in \mathcal{P}(X)$ is $\Pi(\mu, \nu)$. 

The composition of arrows 
is given by composition of transport plans, defined further.  Let
 $\displaystyle \gamma \in \Pi(\mu,
\nu)$ and $\gamma' \in \Pi(\nu, \mu')$. By the disintegration theorem 
\ref{disthm}, there are $\nu$-a.e. uniquely defined measure-valued maps 
$\displaystyle x \in X \mapsto \gamma'_{x}, \gamma_{x} \in \mathcal{P}(X)$ such
that  for any Borel function $\displaystyle 
f: X_{1} \times X_{2} \rightarrow [0,+\infty]$ we have: 
$$\int_{X \times X} f(x_{1}, x_{2}) \mbox{ d} \gamma(x_{1}, x_{2}) \, = \, 
\int_{X} \left( \int_{X} f(x_{1}, x_{2}) \mbox{ d} \gamma_{x_{2}}(x_{1}) \right)
\mbox{ d} \nu(x_{2})$$
$$\int_{X \times X} f(x_{1}, x_{2}) \mbox{ d} \gamma'(x_{1}, x_{2}) \, = \, 
\int_{X} \left( \int_{X} f(x_{1}, x_{2}) \mbox{ d} \gamma'_{x_{1}}(x_{2}) \right)
\mbox{ d} \nu(x_{1})$$
Then the composition $\gamma' \circ
\gamma \in \Pi(\mu, \mu')$ is defined by: 
 for any Borel function $\displaystyle 
f: X_{1} \times X_{2} \rightarrow [0,+\infty]$ 
$$\int_{X \times X} f(x_{1}, x_{3}) \mbox{ d} \left(\gamma' \circ
\gamma \right) (x_{1}, x_{3}) \, = \,  \int_{X} \left( 
\int_{X \times X} f(x_{1}, x_{3}) \mbox{ d} \gamma'_{x_{2}}(x_{3}) \mbox{ d} 
\gamma_{x_{2}}(x_{1}) \right) \mbox{ d} \nu(x_{2})$$
Finally, this category has a contravariant "inverse" functor, which associates 
to each arrow $\displaystyle \gamma \in \Pi(\mu,\nu)$ the arrow 
$\displaystyle \gamma^{-1} \in \Pi(\nu,\mu)$ defined by: for any Borel function $\displaystyle 
f: X_{1} \times X_{2} \rightarrow [0,+\infty]$ 
$$\int_{X \times X} f(x_{1}, x_{2}) \mbox{ d} \gamma^{-1}(x_{1}, x_{2}) \, = \, 
\int_{X \times X} f(x_{1}, x_{2}) \mbox{ d} \gamma(x_{2}, x_{1})$$
\end{defi}

We denote by $\displaystyle Lip_{1}(X)$ the class of 1-Lipschitz maps from $X$
to $\mathbb{R}$. 

\begin{prop}
The category $\displaystyle Trans(X)$ is a separable seminormed category, with
the family $S$ of seminorms: for any $\displaystyle u \in Lip_{1}(X)$ and any 
$\gamma \in Trans(X)$  
$$\rho_{u}(\gamma) \, = \, \mid \int_{X \times X} \left( u(x) - u(y) \right)
\mbox{ d}\gamma(x,y) \mid$$
$\displaystyle Trans(X)$ it is also a normed category, with the norm: 
$$d(\gamma) \, = \, \int_{X \times X} d(x,y) \mbox{ d}\gamma(x,y)$$
\end{prop}

The following relation between the norm $d$ in the transport category and the 
seminorms $\displaystyle \rho_{u}$, is related to Kantorovich formulation 
of the transport problem and the dual of this problem: for any $\displaystyle 
u \in Lip_{1}(X)$ and any $\gamma \in Trans(X)$ we have 
$$d(\gamma) \, \geq \, \rho_{u}(\gamma)$$
The topology of $Trans(X)$ with the norm $d$ is called "strong topology", while 
the topology of $Trans(X)$ as a seminormed groupoid is called the "weak
topology". 

For any pair formed by a Borel function $f: X \rightarrow X$ and by a 
measure $\mu \in \mathcal{P}(X)$, there is an associated transport plan 
$\displaystyle (f, \mu) = (id\, \times f)\sharp \mu \in \Pi(\mu, f\sharp\mu)$. 

\begin{defi}
If $\gamma \in \mathcal{P}(X \times X)$ is representable as a measure 
$(f,\mu)$ then we say that $\gamma$ is induced by a transport map. The category 
$Tmap(X)$ is the subcategory of $Trans(X)$ with objects elements of 
$\mathcal{P}(X)$ and arrows transport plans induced by transport maps. 
\end{defi}

Several facts deserve to be mentioned. 
Notice that if $(f,\mu) = (g,\mu)$ then $f(x) = g(x)$ $\mu$-a.e. in $X$. Also, 
in the category $Tmap(X)$ the composition of arrows is the following: 
$$(g, f\sharp\mu) \circ (f,\mu) \, = \, (g \circ f, \mu)$$
Finally, the sub-category $Invtrans(X)$ of $Tmap(X)$ with objects elements 
of $\mathcal{P}(X)$ and arrows $(f,\mu)$ such that there exists a Borel function
$g$ with $\displaystyle (f,\mu)^{-1} = (g, f\sharp\mu)$, is a groupoid.

\section{Dilations induced deformations of normed groupoids}
Let $(G,d)$ be a  normed groupoid with a separable norm. 
A deformation of $(G,d)$ is basically a "local action" of a commutative 
group $\Gamma$ on $G$ which satisfies several properties.

$(\Gamma, \mid\cdot\mid)$ is a commutative group endowed with a group 
morphism $\mid\cdot\mid : \Gamma \rightarrow (0,+\infty)$ to the 
multiplicative group of positive real numbers. This morphism induces a 
invariant topological filter over $\Gamma$ (a end of $\Gamma$). 
Further we shall write $\varepsilon \rightarrow 0$ for $\varepsilon$ 
converging to this end, and meaning that $\mid\varepsilon\mid\rightarrow 0$. 
The neutral element of $\Gamma$ is denoted by $e$. 

 To any $\varepsilon \in \Gamma$ is associated a transformation 
 $\displaystyle \delta_{\varepsilon}: \, dom(\varepsilon) \rightarrow \, 
 im(\varepsilon)$, which may be  called a dilatation,  dilation, homothety 
  or  contraction. 

For the precise properties of the domains and codomains of $\displaystyle
\delta_{\varepsilon}$ for $\varepsilon \in \Gamma$ see the subsection
\ref{bordet}. For the moment is  sufficient to know that for 
any $\varepsilon \in \Gamma$ we have $\displaystyle Ob(G) = X \subset 
\, dom(\varepsilon)$ and $\displaystyle Ob(G) = X \subset 
\, im(\varepsilon)$. Basically the domain and codomain of $\displaystyle
\delta_{\varepsilon}$ are neighbourhoods of $X$. Moreover, these sets are 
chosen so that various compositions of transformations $\displaystyle
\delta_{\varepsilon}$ are well defined. 

In the formulation of properties of deformations we shall use a uniform
convergence on bounded sets. We explain further what uniform
convergence on bounded sets means in the case of nets of functions indexed with the directed net
the group $\Gamma$ (ordered such that limits are taken in the sense $\mid\varepsilon\mid  
\rightarrow 0$.

\subsection{Introducing dilations}

\begin{defi}
A {\bf dilation of a separated normed groupoid} $(G,d)$ is a map assigning to any 
$\varepsilon \in \Gamma$ a transformation $\delta_{\varepsilon}: \, 
dom(\varepsilon) \rightarrow \, im(\varepsilon)$ which satisfies the following:
\begin{enumerate}
\item[A1.] For any 
$\varepsilon \in \Gamma$ $\alpha \delta_{\varepsilon} = \alpha$. Moreover $\displaystyle
\varepsilon \in \Gamma \mapsto \delta_{\varepsilon}$ is an action of $\Gamma$ on
$G$, that is for any 
$\varepsilon, \mu \in \Gamma$ we have $\displaystyle \delta_{\varepsilon} 
\delta_{\mu} = \delta_{\varepsilon \mu}$, $\displaystyle 
\left(\delta_{\varepsilon} \right)^{-1} = \delta_{\varepsilon^{-1}}$ and 
$\delta_{e} = \, id$. 
\item[A2.] For any $x \in Ob(G)$ and any $\varepsilon \in \Gamma$ we have 
$\displaystyle \delta_{\varepsilon}(x) = x$. Moreover the transformation 
$\displaystyle \delta_{\varepsilon}$ contracts $dom(\varepsilon)$ to $X = Ob(G)$ uniformly on bounded sets, which means that 
the net $\displaystyle d \, \delta_{\varepsilon}$ converges to the constant
function $0$, uniformly on bounded sets, in the sense of definition
\ref{defcon}.
\end{enumerate}
Moreover the domains and codomains $dom(\varepsilon)$,  $im(\varepsilon)$
satisfy the conditions from definition \ref{dax0}, section \ref{bordet}.

A {\bf dilation of a separated normed or seminormed  category} is defined in
the same way as a dilation of a normed groupoid, only that $\alpha$ and
$\omega$ functions are no longer the source and arrow functions, but 
$\alpha(g) = g^{-1}g$ and $\omega(g) = g g^{-1}$. 
\label{ddefor}
\end{defi}

\subsubsection{Deformation of the $\alpha$-double groupoid.} 
The dilation $\delta$  of $(G,d)$ 
induces a right-invariant dilation of the 
normed groupoid $\displaystyle (G\times_{\alpha}D, \tilde{d})$. The proof of the
following proposition is straightforward and we do not write it. 

\begin{prop}
For any $\varepsilon \in \Gamma$ we define  $\displaystyle 
\tilde{\delta}$ on $G\times_{\alpha}G$, given by: 
\begin{equation}
\tilde{\delta}_{\varepsilon}(g,h) =  (\delta_{\varepsilon}(gh^{-1}) h, h)
\label{rtransdel}
\end{equation}
This is a deformation of the normed groupoid  $\displaystyle (G\times_{\alpha}G,
\tilde{d})$ is a normed groupoid and moreover  $dif$
is a morphism of normed groupoids (that is a norm preserving morphism of
groupoids), which commutes with dilations in the sense: for any $\varepsilon
\in \Gamma$ $\displaystyle dif \, \tilde{\delta}_{\varepsilon} \, = \, \delta_{\varepsilon} \, 
dif$. 
\label{pdefindd}
\end{prop} 

\subsubsection{Morphisms of dilations.} Let $\displaystyle (G_{1},d_{1},
\delta^{1})$  and $\displaystyle (G_{2},d_{2},\delta^{2})$ be two dilations 
of the normed groupoids  $\displaystyle (G_{1},d_{1})$, 
$\displaystyle (G_{2},d_{2})$ respectively. Then $\displaystyle F: 
(G_{1},d_{1},\delta^{1}) \rightarrow (G_{2},d_{2},\delta^{2})$ is a morphism of 
dilations if: $F$ is a morphism of groupoids, it preserves the norms (it is 
a isometry) and it commutes with dilations (it is "{\bf linear}").

\subsection{Domains and codomains of dilations}
\label{bordet}

Dilations are locally defined. This is explained in the following definition,
which should be seen as the axiom A0 of dilations.

\begin{defi}
The domains and codomains of a dilation of $(G,d)$ satisfy the following 
Axiom A0: 
\begin{enumerate}
\item[(i)] for any $\varepsilon \in \Gamma$ $\displaystyle Ob(G) = X \subset 
\, dom(\varepsilon)$ and $\displaystyle dom(\varepsilon) = 
dom(\varepsilon)^{-1}$, 
\item[(ii)]  for any bounded set $K \subset Ob(G)$ there are  $1<A<B$  such that 
for any $\varepsilon \in \Gamma$, $\mid \varepsilon\mid \leq 1$: 
$$d^{-1}(\mid \varepsilon\mid) \, \cap \, \alpha^{-1}(K) \, \subset \,  
\delta_{\varepsilon} \, \left( d^{-1}(A) \cap \alpha^{-1}(K) \right)  
\, \subset  \, dom(\varepsilon^{-1}) \, \cap \, \alpha^{-1}(K)  \, \subset $$
\begin{equation}
\subset \,  
\delta_{\varepsilon} \, \left( d^{-1}(B) \cap \alpha^{-1}(K) \right)  \, \subset
\, \delta_{\varepsilon} \, \left( dom(\varepsilon) \cap \alpha^{-1}(K)\right)
\label{eax0}
\end{equation}
\item[(iii)]  for any bounded set $K \subset Ob(G)$ there are $R> 0$ and 
$\displaystyle \varepsilon_{0} \in (0,1]$ such that 
for any $\varepsilon \in \Gamma$, $\displaystyle \mid \varepsilon\mid \leq
\varepsilon_{0}$ and any $\displaystyle g,h \in \mid d^{-1}(R) \cap \,
\alpha^{-1}(K)$  we have:
\begin{equation}
dif(\delta_{\varepsilon} g, \delta_{\varepsilon} h) \, \in \,
dom(\varepsilon^{-1})
\label{edomdif}
\end{equation}
\end{enumerate}
\label{dax0}
\end{defi}

\begin{rk}
Concerning (iii) definition \ref{dax0}, the first part of A1 definition
\ref{ddefor}  
 implies that $\displaystyle dif(\delta_{\varepsilon} g, 
\delta_{\varepsilon}h)$ is well defined for any $\displaystyle 
(g,h) \in G\times_{\alpha}G$ such that $g,h \in dom(\varepsilon)$. 
\end{rk}

\subsection{Induced deformations}

\begin{figure}[h]
\centering
\includegraphics[angle=270, width=0.5\textwidth]{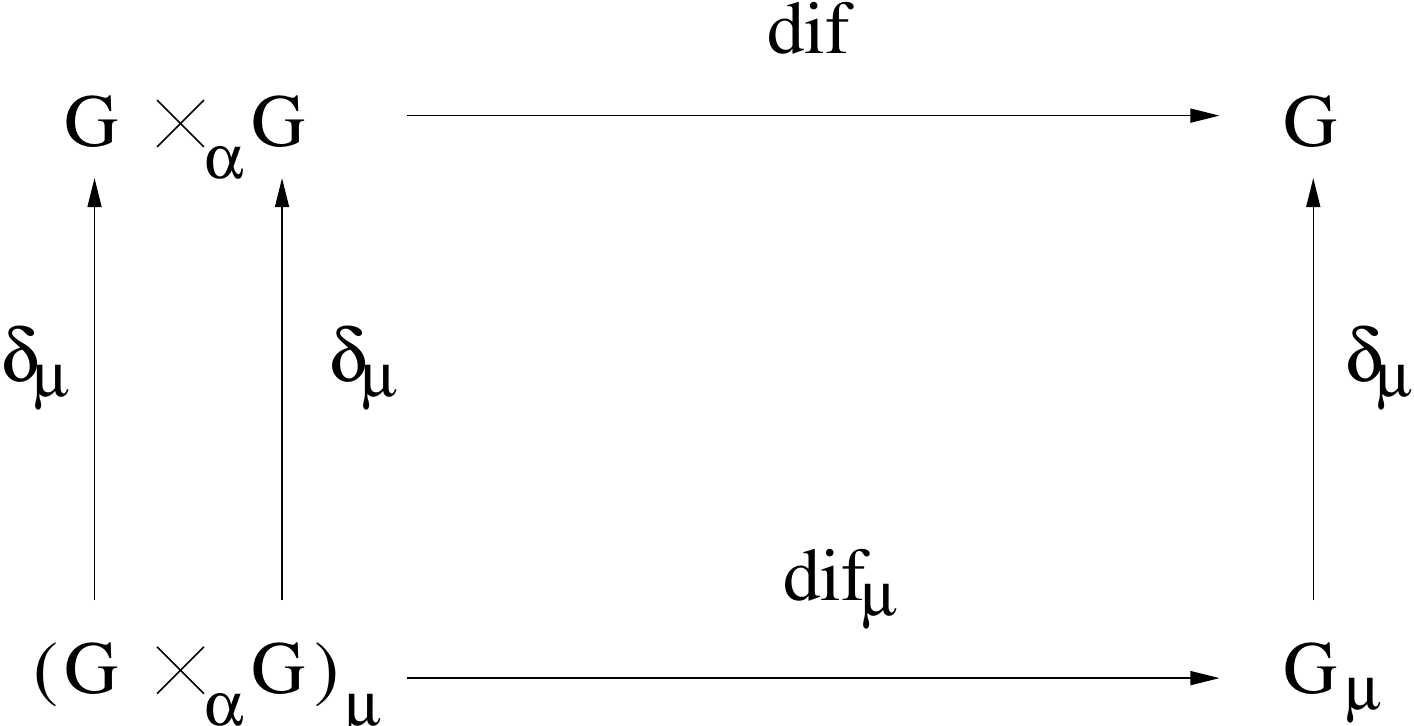}
\caption{Commutative diagram of induced normed groupoids}
\label{figure2}
\end{figure}

The purpose of this section is to define several deformations of normed
groupoids, such that the diagram from figure \ref{figure2} becomes a commutative
diagram of morphisms of dilations. 

Let us consider a triple $(G,d,\delta)$
with $(G,d)$ a normed groupoid and $\delta$  a dilation. For any $\mu \in \Gamma$ there are two normed induced groupoids, such that 
the arrows in the diagram (\ref{figure2}) are morphisms.

\begin{rk}
As dilatations are not globally defined and they are used to transport groupoid 
operations, it follows that the transported objects (operation, norms, ...) 
are not globally defined. Therefore the induced groupoids are not groupoids, 
but "local" groupoids, in a sense which is clear in the context. 
\label{rkdom}
\end{rk}

\begin{defi}
The deformation $\displaystyle (G_{\mu}, d_{\mu}, \delta)$ is equal to $G$ as a set and its 
operations,  norm and dilation are transported by the map 
 $\displaystyle \delta_{\mu} : G_{\mu} \rightarrow G$ (with the precautions 
 concerning the domains of definition of the transported objects mentioned in 
 remark \ref{rkdom}). 
 
 The deformation $\displaystyle \left(\left(G\times_{\alpha}G\right)_{\mu},
 \tilde{d}_{\mu}, \tilde{\delta}_{\mu}\right)$ is equal to 
 $\displaystyle G\times_{\alpha}G$ as a groupoid and its norm and dilation 
 are transported by the map 
 $\displaystyle \delta_{\mu} \times \delta_{\mu} : \left(G\times_{\alpha}G\right)_{\mu} 
 \rightarrow G\times_{\alpha}G$.
 \label{defdia}
 \end{defi}
 More precisely, the deformation $\displaystyle (G_{\mu}, d_{\mu}, \delta)$  is described by: 
 \begin{enumerate}
 \item[-] $\displaystyle G_{\mu} = G$  as a set, $\displaystyle 
 \alpha_{\mu} = \alpha$ and $\displaystyle \omega_{\mu} = \omega \,
 \delta_{\mu}$, which follow from the computations using A1, A2 definition
 \ref{ddefor}: 
 $$\alpha_{\mu} = \delta_{\mu}^{-1} \, \alpha \, \delta_{\mu} =
 \delta_{\mu}^{-1} \, \alpha = \alpha \quad , \quad 
 \omega_{\mu} = \delta_{\mu}^{-1} \, \omega \, \delta_{\mu} =
 \omega \, \delta_{\mu}$$
 Also $\displaystyle Ob(G_{\mu}) = Ob(G)$. 
 \item[-] the composition operation and inverse are 
 $$m_{\mu}(g,h) \, = \, \delta_{\mu}^{-1} \left( \delta_{\mu}(g) \,
 \delta_{\mu}(h)\right) \quad , \quad inv_{\mu}\, g \, = \, \delta_{\mu}^{-1} \,
 inv \, \delta_{\mu} \, (g)$$
 These are well defined (at least locally) because of the axiom A0 definition
 \ref{dax0}. Notice that $\displaystyle dif_{\mu}$ from the diagram 
 \ref{figure2} appears as   the difference function 
 associated with the operation $\displaystyle m_{\mu}$, defined as 
\begin{equation}
 dif_{\mu} : G_{\mu} \times_{\alpha_{\mu}} G_{\mu} \rightarrow G_{\mu} \quad , 
 \quad dif_{\mu}(g,h) \, = \, \delta_{\mu}^{-1} \left( \delta_{\mu}(g) \,
 \left(\delta_{\mu}(h)\right)^{-1} \right)
 \label{dfird}
 \end{equation}
\item[-] the norm $\displaystyle d_{\mu}$ is defined as: 
\begin{equation}
d_{\mu}(g) \, = \, \frac{1}{\mid\mu\mid} \, d(\delta_{\mu} g)
\label{dfirdem}
\end{equation}
\item[-] we may transport the dilation $\delta$ of $(G,d)$ into a dilation
$\displaystyle \delta_{\mu, \cdot}$ of $\displaystyle (G_{\mu}, d_{\mu})$, but
from the commutativity of $\Gamma$ we get that 
$$\delta_{\mu, \varepsilon} \, = \, \delta_{\mu}^{-1} \, \delta_{\varepsilon} \,
\delta_{\mu} \, = \, \delta_{\varepsilon}$$
therefore it is the same dilation.
\end{enumerate}

The deformation $\displaystyle \left(\left(G\times_{\alpha}G\right)_{\mu},
 \tilde{d}_{\mu}, \tilde{\delta}_{\mu}\right)$   is described by: 
 \begin{enumerate}
 \item[-] $\displaystyle \left(G\times_{\alpha}G\right)_{\mu} =
 G\times_{\alpha}G$  as a groupoid; remark that this is compatible with the
 transport of operations using the map  $\displaystyle \delta_{\mu} \times
 \delta_{\mu}$ (because this map is an endomorphism of the groupoid 
 $\displaystyle G\times_{\alpha}G$), 

 \item[-] with respect to the relation (\ref{dfird}), notice that 
 $\displaystyle G_{\mu} \times_{\alpha_{\mu}} G_{\mu} = G\times_{\alpha}G = 
 \left(G\times_{\alpha}G\right)_{\mu}$ and $\displaystyle dif_{\mu}$ as
 represented in figure \ref{figure2} is a 
 morphism of groupoids, 
 
\item[-] the norm $\displaystyle \tilde{d}_{\mu}$ is defined as: 
\begin{equation}
\tilde{d}_{\mu}(g,h) \, = \, \frac{1}{\mid\mu\mid} \, \tilde{d}(\delta_{\mu} g, 
\delta_{\mu} h) 
\label{dsecdem}
\end{equation}
and it is easy to check that $\displaystyle dif_{\mu}$ is also a isometry. 
\item[-] we  transport the dilation $\displaystyle\tilde{\delta}$ of 
$\displaystyle \left(G\times_{\alpha}G, \tilde{d}\right))$ into a dilation
$\displaystyle \tilde{\delta}_{\mu, \cdot}$ 
\begin{equation}
\tilde{\delta}_{\mu, \varepsilon}(g,h) \, = \, \left( \delta_{\mu^{-1}} \left( 
\delta_{\varepsilon} \left( \delta_{\mu}(g) \, \left( \delta_{\mu} (h)
\right)^{-1} \right) \delta_{\mu}(h) \right) \, , \, h \right)
\label{dsecdef}
\end{equation}
\end{enumerate} 
The commutativity of the diagram \ref{figure2} is clear now.

\section{Algebraic operations from dilations}

At the core of the introduction of dilations lies the fact that 
we can construct group operations from them. More precisely we are able to 
construct, by using compositions of dilations and the groupoid operation, 
approximately associative operations which shall lead us eventually to group 
operations in the tangent groupoid of a dilation.

\subsection{A general construction}

Let $(G,d,\delta)$ be a dilation and $\displaystyle \left( G\times_{\alpha}G
, \tilde{d}, \tilde{\delta}\right)$ the associated dilation of the
$\alpha$-double groupoid. Further we shall be interested only in the properties
of the following map.

\begin{defi}
 For any $x \in Ob(G)$ and any $\varepsilon \in \Gamma$ 
we define the {\bf dilatation}: 
\begin{equation}
\delta_{\varepsilon}^{(\cdot )} \, (\cdot) : \alpha^{-1}(x) \times \alpha^{-1}(x)
 \rightarrow \alpha^{-1}(x) \quad , \quad 
\delta^{h}_{\varepsilon} g \, = \, \delta_{\varepsilon} \left(g \, h^{-1}\right)
h
\label{defdila}
\end{equation}
\label{defdilc}
\end{defi}

\begin{rk}
The domain of definition of $\displaystyle \delta_{\varepsilon}^{(\cdot )} \,
( \cdot )$ is in fact only a subset of $\displaystyle \alpha^{-1}(x) \times 
\alpha^{-1}(x)$, according to the Axiom 0 explained in definition \ref{dax0} 
section \ref{bordet}. 
\end{rk}

This map comes from the definition (\ref{rtransdel}) of the dilation 
$\displaystyle \tilde{\delta}$, namely 
\begin{equation}
\tilde{\delta}_{\varepsilon}(g,h) \, = \, \left( \delta^{h}_{\varepsilon} g , 
h \right)
\label{defdilb}
\end{equation}

For any $\varepsilon \in \Gamma$ with $\mid\varepsilon\mid$ sufficiently small
we can define $\displaystyle dif_{\varepsilon}$ (as in figure \ref{figure2})
from (a subset of)  $\displaystyle 
 G\times_{\alpha}G$ to $G$. 
Remark that $\displaystyle \alpha(dif_{\varepsilon}(g,h)) =
\omega(\delta_{\varepsilon} h)$, therefore the following  composition is well 
defined:
\begin{equation}
\Delta_{\varepsilon}(g,h) = dif_{\varepsilon}(g,h) \, \delta_{\varepsilon}h
\label{difdef}
\end{equation}
Then $\displaystyle \alpha \, \Delta_{\varepsilon}(g,h) = \alpha (g) = 
\alpha(h)$. 

Related to the function $\displaystyle \Delta_{\varepsilon}$ is the following 
\begin{equation}
inv_{\varepsilon}(g) \, = \, \Delta_{\varepsilon}(\alpha(g), g) 
\label{invdef}
\end{equation}  

The following expression makes sense too, for any pair of elements $(g,h)$ from 
(a subset of) $\displaystyle  G\times_{\alpha}G$: 
\begin{equation}
\Sigma_{\varepsilon}(g,h) = \delta_{\varepsilon^{-1}} \left[
\delta_{\varepsilon} \left( g \left( \delta_{\varepsilon} h \right)^{-1} 
\right) \delta_{\varepsilon} h \right] 
\label{sumdef}
\end{equation}
It is also true that $\displaystyle \alpha \, \Sigma_{\varepsilon}(g,h) = \alpha (g) = 
\alpha(h)$.

These three functions are interesting operations. The function $\displaystyle 
\Delta_{\varepsilon}$ is an approximate difference operation, $\displaystyle 
inv_{\varepsilon}$ is an approximate inverse and $\displaystyle 
\Sigma_{\varepsilon}$ is an approximate sum operation.

\subsubsection{A graphic construction of approximate difference operation}
\begin{figure}[h]
\centering
\includegraphics[angle=270, width=0.5\textwidth]{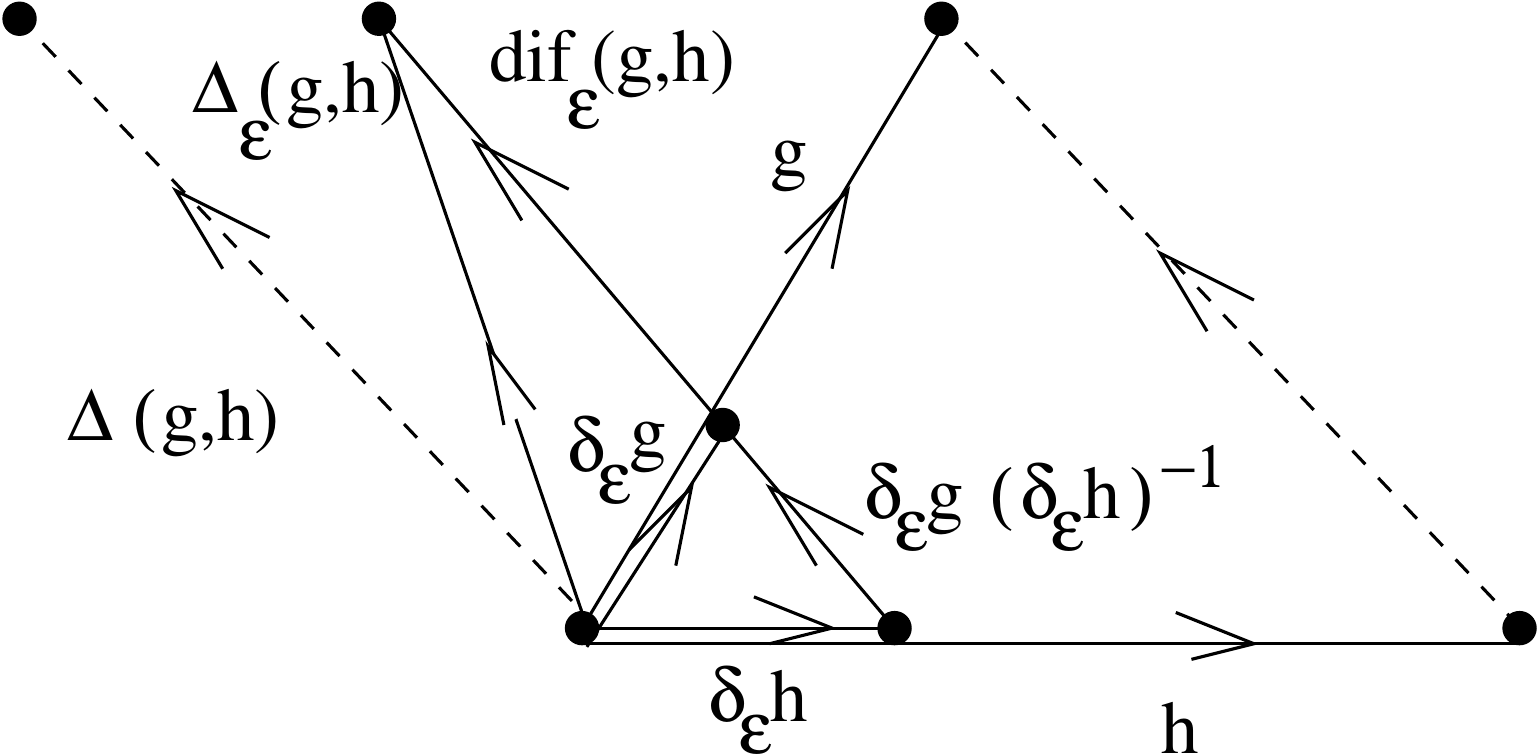}
\caption{The meaning of the "approximate difference" operation $\displaystyle 
\Delta_{\varepsilon}(g,h)$.}
\label{figure1}
\end{figure}

A look at the figure 
\ref{figure1} will help. There is graphically explained how 
$\displaystyle \Delta_{\varepsilon}(g,h)$ is constructed. 

Let us imagine 
that we are looking at a figure in the Euclidean plane. Then $\displaystyle 
\delta_{\varepsilon}$ is just a homothety, $g$, $h$ are vectors with the 
same origin $\alpha(g) = \alpha(h)$, $\displaystyle \Delta(g,h)$ is the 
difference of vectors $-g +h$ (or $h -g$, it's the same as long as we are in a
commutative world). In the Euclidean plane, as $\mid\varepsilon\mid$ goes to 
$0$, the "vector" $\displaystyle dif_{\varepsilon}(g,h)$ slides towards $\Delta(g,h)$ 
 and $\displaystyle \Delta_{\varepsilon}(g,h)$ is obtained from 
$\displaystyle dif_{\varepsilon}(g,h)$ by composition with the vector 
$\displaystyle \delta_{\varepsilon}h$. Thus $\displaystyle 
\Delta_{\varepsilon}(g,h)$ has the meaning of a approximate difference of 
vectors $g$, $h$.

\subsection{Idempotent right quasigroup and induced operations}

\subsubsection{Approximate operations from dilatations} 
The functions $\displaystyle 
\Delta_{\varepsilon}$, $\displaystyle 
inv_{\varepsilon}$  and $\displaystyle 
\Sigma_{\varepsilon}$ can be expressed in terms of dilatations introduced in 
definition \ref{defdilc}. Indeed, let us define, for any triple $u,g,h \in G$ 
with $\alpha(u)= \alpha(g) = \alpha(h)$, and such that $d(u), d(g), d(h)$ are
sufficiently small, the following approximate difference function with three
arguments: 
\begin{equation}
\Delta_{\varepsilon}^{u}(g,h) \,  =  \,
\delta_{\varepsilon^{-1}}^{\delta_{\varepsilon}^{u} g}
\delta_{\varepsilon}^{u} h
\label{def3dif}
\end{equation}
the approximate inverse function with two arguments: 
\begin{equation}
inv_{\varepsilon}^{u}(g) \,  =  \,
\delta_{\varepsilon^{-1}}^{\delta_{\varepsilon}^{u} g} u \, = \,
\Delta_{\varepsilon}^{u}(g,u)
\label{def3inv}
\end{equation}
and the following approximate sum function with three arguments: 
\begin{equation}
\Sigma_{\varepsilon}^{u}(g,h) \,  =  \,
\delta_{\varepsilon^{-1}}^{u} \, \delta_{\varepsilon}^{\delta_{\varepsilon}^{u}
g} h
\label{def3sum}
\end{equation}
We have then: 
\begin{equation}
\Delta_{\varepsilon}(h u^{-1},g u^{-1}) \, = \, 
\Delta^{u}_{\varepsilon} (g,h) u^{-1}) \quad ,
\quad \Sigma_{\varepsilon}(h u^{-1},g u^{-1}) \, = \, 
\Delta^{u}_{\varepsilon} (g,h) u^{-1})
\label{compa}
\end{equation}

\subsubsection{$\Gamma$-idempotent right quasigroups} We are in the framework of emergent algebras and idempotent right quasigroups,
as introduced in \cite{buligairq}. We recall here the definition of a 
idempotent right quasigroup and induced operations.

\begin{defi} An idempotent right quasigroup (irq) is a  set $X$ endowed with two  operations $\circ$ and $\bullet$, which satisfy the following axioms: for any $x , y \in X$  
\begin{enumerate}
\item[(P1)] \hspace{2.cm} $\displaystyle x \, \circ \, \left( x\, \bullet \,  y \right) \, = \, x \, \bullet \, \left( x\, \circ \,  y \right) \, = \, y$
\item[(P2)] \hspace{2.cm} $\displaystyle x \, \circ \, x \, = \, x \, \bullet \, x \,  = \,  x$
\end{enumerate}
We use these operations to define the sum, difference and inverse operations of the irq:  for 
any $x,u,v \in X$ 
\begin{enumerate}
\item[(a)] the difference operation is $\displaystyle (xuv) \, = \, \left( x \,
\circ \, u \right) \, \bullet \, \left( x \, \circ \, v \right)$. 
 By fixing the first variable $x$ we obtain the difference operation based at $x$: 
$\displaystyle v \, -^{x} \, u \, =  \, dif^{x}(u,v) \, = \, (xuv)$. 
\item[(b)] the sum operation is 
$\displaystyle )xuv( \, = \, x \, \bullet \left( \left( x \, \circ \, u \right) \, \circ \, 
 v \right)$. 
  By fixing the first variable $x$ we obtain the sum operation based at $x$: 
$\displaystyle u \, +^{x} \, v \, =  \, sum^{x}(u,v) \, = \, )xuv($.  
\item[(a)] the inverse operation is  
$\displaystyle inv(x,u) \, = \, \left( x \, \circ \, u \right) \, \bullet \,  x
$. By fixing the first variable $x$ we obtain the inverse  operator based 
at $x$: $\displaystyle  -^{x} \, u \, =  \, inv^{x} u \, = \, inv(x,u)$. 
\end{enumerate}
 For any $\displaystyle 
k \in \mathbb{Z}^{*} = \mathbb{Z} \setminus \left\{ 0 \right\}$  we define also 
the following operations: 
\begin{enumerate}
\item[-] $\displaystyle x \, \circ_{1}\,  u \, = \, x \, \circ \, u$,  $\displaystyle x \, \bullet_{1}\,  u \, = \, x \, \bullet \, u$, 
\item[-]  for any $k > 0$ let 
$\displaystyle  x \, \circ_{k+1}\,  u \, = \, x \, \circ \left(x \circ_{k} \, u \right)$ and 
$\displaystyle  x \, \bullet_{k+1}\,  u \, = \, x \, \bullet \left(x \bullet_{k} \, u \right)$, 
\item[-] for any $k < 0$ let 
$\displaystyle  x \, \circ_{k}\,  u \, = \, x \bullet_{-k} \, u$ and 
$\displaystyle  x \, \bullet_{k}\,  u \, = \, x \circ_{-k} \, u $. 
\end{enumerate}

\label{dplay}
\end{defi}

For any $\displaystyle k \in \mathbb{Z}^{*}$  the triple 
$\displaystyle (X, \circ_{k} , \bullet_{k} )$ is a irq. We denote  the difference, sum and inverse operations of $\displaystyle (X, \circ_{k} , \bullet_{k} )$  by the same symbols as 
the ones used for  $(X, \circ , \bullet )$,  with  a subscript "$k$". 

For any $\varepsilon \in \Gamma$  and for any $x \in X$ 
we can define a irq operation on $\displaystyle \alpha^{-1}(x)$ by 
$\displaystyle g \, \circ_{\varepsilon} \, h \, = \, \delta_{\varepsilon}^{g}
h$. We have then: 
$$u \, +^{g} \, h  \, = \, \Sigma^{g}_{\varepsilon}(u , h) \quad , \quad 
u \, -^{g} \, h  \, = \, \Sigma^{g}_{\varepsilon}(h , u)$$

By computation it follows that $\displaystyle
\left(\circ_{\varepsilon}\right)_{k} \, = \, \circ_{\varepsilon^{k}}$. The
approximate difference, sum and inverse operations are exactly the ones
introduced in the preceding section. 
In \cite{buligairq} we introduced idempotent right quasigroups and then iterates
of the operations indexed by a parameter $\displaystyle k \in \mathbb{N}$. This
was done in order to simplify the notations mostly. Here, in the presence of the
group $\Gamma$, we might define a $\Gamma$-irq. 

\begin{defi}
Let $\Gamma$ be a commutative group. A $\Gamma$-idempotent right quasigroup 
is a set $X$ with a function $\displaystyle \varepsilon \in \Gamma \mapsto 
\circ_{\varepsilon}$ such that $\displaystyle (X, \circ_{\varepsilon})$ is a irq
and moreover for any $\varepsilon, \mu \in \Gamma$ and any $x, y \in X$ we have 
$$x \, \circ_{\varepsilon} \, \left( x \, \circ_{\mu} \, y \right) \, = \, 
x \circ_{\varepsilon \mu}$$
\label{defgammairq}
\end{defi}

It is then obvious that if $(X, \circ)$ is a irq then $\displaystyle 
(X, k \in \mathbb{Z} \mapsto \circ_{k})$ is a $\mathbb{Z}$-irq (we define 
$\displaystyle x \, \circ_{0} \, y \, = \, y$).

The following is a slight modification of proposition 3.4 and point (k)
proposition 3.5 \cite{buligairq}, for
the case of $\Gamma$-irqs (the proof of this proposition is almost 
identical, with obvious modifications, with the proof of the original
proposition).

\begin{prop}
In any irq  $\displaystyle (X,\circ_{\varepsilon})_{\varepsilon \in \Gamma}$ be
a $\Gamma$-irq. Then we have the relations: 
\begin{enumerate}
\item[(a)] $\displaystyle \left( u \, +^{x}_{\varepsilon} \, v \right) \,
-^{x}_{\varepsilon} \, u \, = \, v $
\item[(b)] $\displaystyle u \, +^{x}_{\varepsilon} \, \left( v \, -^{x}_{\varepsilon} \, u \right) \, = \, v $
\item[(c)] $\displaystyle v \, -^{x}_{\varepsilon} \, u \, = \, \left(-^{x}_{\varepsilon} u\right)  \, +^{x \circ u}_{\varepsilon} \, v  $
\item[(d)] $\displaystyle -^{x\circ u}_{\varepsilon} \, \left( -^{x}_{\varepsilon} \, u \right) \, = \, u $
\item[(e)] $\displaystyle u \, +^{x}_{\varepsilon} \, \left( v \, +^{x\circ u}_{\varepsilon} \, w \right) \, = \, \left( u \, +^{x}_{\varepsilon} \, v \right) \, +^{x}_{\varepsilon} \, w $
\item[(f)] $\displaystyle  -^{x}_{\varepsilon} \, u \, = \,  x \, -^{x}_{\varepsilon} \, u $
\item[(g)] $\displaystyle  x \, +^{x}_{\varepsilon} \, u \, = \,  u $
\item[(k)] for any $\displaystyle \varepsilon, \mu \in \mathbb{Z}^{*}$ and any $x, u , v \in X$ we have the distributivity property: 
$$\displaystyle (x \circ_{\mu} v) \, -_{\varepsilon}^{x} \, ( x \circ_{\mu} u)
\, = \, \left( x \circ_{\varepsilon \mu} u \right) \, 
\circ_{\mu} \, \left( v \, -_{\varepsilon \mu}^{x} \, u \right)$$
\end{enumerate}
\label{pplay}
\end{prop}

Later we shall apply  this proposition for the irq $\displaystyle
\alpha^{-1}(x)$ with the operations induced by dilatations $\displaystyle
\delta_{\varepsilon}$.

\section{Limits of induced deformations}

As $\displaystyle \mid\mu\mid
\rightarrow 0$ the components of  the deformations indexed by 
$\mu$ from the diagram \ref{figure2} 
(namely the operation, norm and respective dilation maps) may converge in the 
sense of section \ref{secdefconv} to the components of  another
normed groupoid with dilations.

\subsection{The weak case: dilatation structures on metric spaces}

This is the case when only $\displaystyle \left(\left(G\times_{\alpha}G\right)_{\mu},
 \tilde{d}_{\mu}, \tilde{\delta}_{\mu}\right)$ converges. There is no condition 
 of convergence upon $\displaystyle (G_{\mu}, d_{\mu}, \delta)$, nor upon 
 the difference function $\displaystyle dif_{\mu}$.

\begin{defi}
A dilation $(G,d,\delta)$ is a {\bf groupoid weak $\delta$-structure} (
gw $\delta$-structure) 
if it satisfies the following two axioms: 
\begin{enumerate}
\item[A3.] There is a  function   
$\displaystyle \bar{d}: G\times_{\alpha}G \cap U^{2} \rightarrow \mathbb{R}$ 
which is the limit  
$$ \lim_{\varepsilon \rightarrow 0} \frac{1}{\mid\varepsilon\mid} \, 
d \, dif(\delta_{\varepsilon} g, \delta_{\varepsilon} h) \, = \, \tilde{d}_{0}(g,h)$$
uniformly on bounded sets in the sense of 
definition \ref{defcon}. Moreover the convergence with respect to 
$\bar{d}$ is the same as the convergence with respect to $\tilde{d}$ and in
particular $\displaystyle \tilde{d}_{0}(g,h) = 0$ implies
$g = h$. 
\item[A4weak.] There is a dilation $\displaystyle \bar{\delta}$ of the normed
groupoid $\displaystyle (G\times_{\alpha}G, \tilde{d}_{0})$ such that 
for any $\varepsilon \in \Gamma$ the transformation 
$\displaystyle \tilde{\delta}_{\mu, \varepsilon}$ converges uniformly on bounded
sets to $\displaystyle \bar{\delta}_{\varepsilon}$.
 \end{enumerate}
\label{defgdsweak}
\end{defi}

\begin{rk}
For a gw $\delta$-structure the function $\bar{d}$ has the following 
properties of a distance: for any $\displaystyle (g,h) \in G\times_{\alpha}G \cap U^{2}$
\begin{enumerate}
\item[(a)] $\bar{d}(g,h) = 0$ if and only if $g = h$, 
\item[(b)]  $\displaystyle \bar{d}(g,h) \leq \bar{d}(g) + \bar{d}(h)$, 
\item[(c)] $\bar{d}(g,h) = \bar{d}(h,g)$.
\end{enumerate}
This means that for any $x \in Ob(G)$ the function $\bar{d}$ gives 
a distance on the set $\displaystyle \alpha^{-1}(x) \cap U$. 
\label{rema2}
\end{rk}
Indeed, these properties of the function $\bar{d}$ come from the following 
observation. Let us define on $\displaystyle G\times_{\alpha}G$ the function:
$$d(g,h) = d(gh^{-1}) = d \, dif(g,h)$$
Then for any $x \in Ob(G)$ the function $d$ (with two arguments) gives 
a distance on the set $\displaystyle \alpha^{-1}(x)$. 
In the case of a $\delta$-structure the axiom A3 can be written as: 
\begin{equation}
\lim_{\varepsilon \rightarrow 0} \frac{1}{\mid\varepsilon\mid} \, 
d(\delta_{\varepsilon} g, \delta_{\varepsilon} h) \, = \, \bar{d}(g,h)
\label{a3inter}
\end{equation}
uniformly on bounded sets. This gives properties  (b), (c) above from 
 a passage to the limit of the properties of the distance $d$.

For any $x \in Ob(G)$ the restriction of the norm $\tilde{d}$ on the trivial groupoid 
$\displaystyle \alpha^{-1}(x)\times\alpha^{-1}(x)$ gives a distance on 
the space $\displaystyle \alpha^{-1}(x)$. The dilatation
$\tilde{\delta}$ has the property: for any $\varepsilon \in \Gamma$ and 
$x \in Ob(G)$ 
$$\tilde{\delta}_{\varepsilon} \, \alpha^{-1}(x) \subset \alpha^{-1}(x)$$
therefore we can define $\displaystyle \delta^{h}_{\varepsilon}$ from 
(a subset of) $\displaystyle \alpha^{-1} ( \alpha(h))$ to 
$\displaystyle \alpha^{-1} ( \alpha(h))$ by: 
\begin{equation}
\delta^{h}_{\varepsilon} g = \delta_{\varepsilon}(gh^{-1}) h 
\label{firstcomright}
\end{equation}

\begin{thm}
Suppose that $(G,d,\delta)$ is a gw $\delta$-structure. Then 
for any $x \in Ob(G)$ the triple $\displaystyle (\alpha^{-1}(x), \tilde{d},
\delta)$ is a  dilatation structure, with $\delta$ defined by (\ref{firstcomright}) and $\displaystyle
\tilde{d}$ restrictioned to $\displaystyle \alpha^{-1}(x)$. 
\end{thm}

The proof is just a translation of the definition \ref{defgdsweak} in terms of 
metric spaces, using the equivalence between metric spaces and normed trivial
groupoids. At the end we obtain definition \ref{defweakstrong} 
of dilatation structures on metric spaces, given further.

\subsubsection{Dilatation structures on metric spaces.} 
For simplicity we shall list the   axioms of  a dilatation structure $(X,d,\delta)$ without concerning about
domains and codomains of dilatations. For the full definition of dilatation
structure, as well as for their main properties and examples, see
\cite{buligadil1}, \cite{buligadil2}, \cite{buligasr}. The notion appeared from
my efforts to understand  the last section of the paper  
\cite{bell} (see also \cite{pansu}, \cite{gromovsr}, \cite{marmos1},
\cite{marmos2}). 

However, notice several differences with respect to the original definition of 
dilatation structures: 
\begin{enumerate}
 \item[(a)] in the following definition \ref{defweakstrong} we are no longer 
 asking the metric space $(X,d)$ to be locally compact. 
 Also, uniform convergence in compact sets is replaced by uniform convergence in bounded sets. 
\item[(b)] because of the modifications explained at (a), we have to ask
explicitly that the uniformities induced by $\displaystyle d^{x}$ and $d$ 
are the same. 
\item[(c)] finally, dilatation structures in the sense of the following
definition \ref{defweakstrong} are a bit stronger than dilatation structures in
the sense introduced and studied in \cite{buligadil1}, \cite{buligadil2}, namely
we ask for the existence of a "limit dilatation", see the last axiom. This limit
exists for strong dilatation structures, but not for dilatation structures 
in the sense introduced in \cite{buligadil1}, \cite{buligadil2}. 
\end{enumerate}

\begin{defi}
A triple $(X, d, \delta)$ is a {\bf dilatation structure} if $(X, d)$ is a
 metric space and 
 $$\delta: \Gamma \times \left\{ (x,y) \in X \times X \mbox{ : } y \in
 dom(\varepsilon, x) \right\} \rightarrow X \quad , \quad \delta(\varepsilon, x,
 y) \, = \, \delta^{x}_{\varepsilon} y $$
 is a function with the following properties: 
 
\begin{enumerate}
\item[{\bf A1.}]  For any point $x \in X$ the function $\delta$ induces an
action  $\displaystyle \delta^{x}: \Gamma \rightarrow
End(X,d, x)$, where $End(X,d, x)$ is the collection of all continuous, with continuous inverse transformations
$\phi: (X,d) \rightarrow (X,d)$ such that $\phi(x) = x$.

\item[{\bf A2.}]
 The function $\displaystyle \delta $ is continuous. Moreover, it can be
continuously extended to $\displaystyle \bar{\Gamma} \times X \times X$ by $\delta (0,x, y) = x$  and the limit
$$\lim_{\varepsilon\rightarrow 0} \delta_{\varepsilon}^{x} y \, = \, x  $$
is uniform with respect to $x,y$ in bounded set.

\item[{\bf A3.}] There is $A > 1$ such that  for any $x$ there exists
 a  function $\displaystyle (u,v) \mapsto d^{x}(u,v)$, defined for any
$u,v$ in the closed ball (in distance d) $\displaystyle
\bar{B}(x,A)$, such that
$$\lim_{\varepsilon \rightarrow 0} \quad \sup  \left\{  \mid \frac{1}{\mid 
\varepsilon \mid} d(\delta^{x}_{\varepsilon} u,
\delta^{x}_{\varepsilon} v) \ - \ d^{x}(u,v) \mid \mbox{ :  } u,v \in \bar{B}_{d}(x,A)\right\} \ =  \ 0$$
uniformly with respect to $x$ in bounded set. 
Moreover the uniformity induced by $d^{x}$ is the same as the uniformity induced
by $d$, in particular $\displaystyle d^{x}(u,v) = 0$ implies $u = v$. 
\item[{\bf A4weak.}] (for metric spaces) The following limit exists: 
$$\lim_{\varepsilon \rightarrow 0} \delta^{x}_{\varepsilon^{-1}} \, 
\delta^{\delta^{x}_{\varepsilon} u}_{\mu} \, \delta^{x}_{\varepsilon} v \, = \, 
\bar{\delta}^{x, u}_{\mu} v $$
for any $\mu \in \Gamma$, uniformly with respect to $x, u , v$ in bounded sets. 

\end{enumerate}

\label{defweakstrong}
\end{defi}

\begin{rk}
In particular the axiom A1 tells us that $\displaystyle \delta^{x}_{\varepsilon} x = x$ for any $x \in X$, $\varepsilon \in \Gamma$,
also $\displaystyle \delta^{x}_{1} y = y$ for any $x,y \in X$, and $\displaystyle \delta^{x}_{\varepsilon} \delta^{x}_{\mu}
 y = \delta^{x}_{\varepsilon \mu} y$ for any $x, y \in X$ and $\varepsilon, \mu \in \Gamma$.
\end{rk}
\begin{rk}
In axiom A2 we may alternatively put that the limit is uniform with respect to 
$d(x,y)$. Similarly, we may ask in axiom A4weak (for metric spaces) that the
limit is uniform with respect to $d(x, u)$, $d(x, v)$. 
\end{rk}

\begin{rk}
It is easy to see that:
\begin{enumerate}
\item[(a)] If $(X,d)$ is locally compact then the function  $\displaystyle
d^{x}$ is continuous as an uniform limit of continuous functions on a compact
set. If $(X,d)$ is also separable then from the existence of the limit
$\displaystyle d^{x}$ and from axiom A1 we obtain the fact that $\displaystyle
d^{x}$ and $d$ induce the same uniformities. 
\item[(b)] By definition $\displaystyle d^{x}$ is symmetric 
and  satisfies the triangle inequality, but it can be a degenerated distance 
function: there might exist  $\displaystyle v,w $ such that $\displaystyle
d^{x}(v,w) = 0$.But the end of axiom A2 eliminates this possibility.  
\end{enumerate}
\end{rk}

\begin{prop}
Let $(X,d, \delta)$ be a dilatation structure,  $x \in X$, and let 
$$\delta^{x}_{\varepsilon} \, d (u, v) \, = \, \frac{1}{\mid \varepsilon \mid}
\, d( \delta^{x}_{\varepsilon} u , \delta^{x}_{\varepsilon} v )$$
 Then the net of metric spaces $\displaystyle (\bar{B}_{d}(x,A),
 \delta^{x}_{\varepsilon} d)$ converges in the Gromov-Hausdorff sense to the 
 metric space $\displaystyle (\bar{B}_{d}(x,A), d^{x})$. Moreover this metric
 space is a metric cone, in the following sense: for any $\mu \in \Gamma$ such
 that $\mid \mu \mid < 1$ we have $\displaystyle \bar{\delta}^{x, x}_{\mu} \, =
 \, \delta^{x}_{\mu}$ and 
 $$d^{x} ( \delta^{x}_{\mu} u , \delta^{x}_{\mu} v ) \, = \, \mid \mu \mid 
 \, d^{x}(u,v)$$
 \end{prop}
 
 \paragraph{Proof.} 
 The first part of the proposition is just a reformulation  of axiom A3, 
 without the condition of uniform convergence. For the second part 
 remark  that 
 $$\delta^{x}_{\varepsilon^{-1}} \, 
\delta^{\delta^{x}_{\varepsilon} x}_{\mu} \, \delta^{x}_{\varepsilon} v \, = \, 
\delta_{\mu}^{x} v$$
and also that  
 $$\frac{1}{\mid \varepsilon \mid} d( \delta^{x}_{\varepsilon} \, 
 \delta^{x}_{\mu} u , \delta^{x}_{\varepsilon} \, 
 \delta^{x}_{\mu} v ) \, = \, \mid \mu \mid \, \delta^{x}_{\varepsilon \mu} \, d
 (u,v) $$
 Therefore if we pass to the limit 
 with $\varepsilon \rightarrow 0$ in these two relations we get the desired
 conclusion. \quad $\square$
 
\subsubsection{The translation groupoid associated to a dilatation structure}

We may associate to a dilatation structure a groupoid which is not normed, which
look similar to a Riemannian groupoid (section \ref{riemgrou}). 

\begin{defi}
Let $(X,d,\delta)$ be a dilatation structure. The  translation groupoid 
$Tr(X,d,\delta)$ has as objects the distances 
$$d^{x}_{\varepsilon} (u,v) \, \, \frac{1}{\mid \varepsilon \mid} d\left(
\delta^{x}_{\varepsilon} u , \delta^{x}_{\varepsilon} v \right)$$ 
for all $x \in X$ and all $\varepsilon \in \Gamma$. The arrows are of the form 
$\displaystyle \Sigma_{\varepsilon}^{x}(u, \cdot)$, with 
$$\alpha \left( \Sigma_{\varepsilon}^{x}(u, \cdot) \right) \, = \,
d^{\delta^{x}_{\varepsilon} u}_{\varepsilon} \quad , \quad 
\omega \left( \Sigma_{\varepsilon}^{x}(u, \cdot) \right) \, = \,
d^{x}_{\varepsilon}$$
Composition of arrows is composition of functions. 
\label{transgrou}
\end{defi}

\begin{prop}
$Tr(X,d,\delta)$ is a groupoid. Moreover, arrows are isometries, in the sense 
that for any $\varepsilon \in \Gamma$ and $x, u \in X$ the arrow 
$\displaystyle \Sigma^{x}_{\varepsilon}(u, \cdot )$ is an isometry from the
source to the range: for ant $v, w \in X$ we have: 
$$d^{x}_{\varepsilon}\left( \Sigma^{x}_{\varepsilon}(u,v) ,
\Sigma^{x}_{\varepsilon}(u,w) \right) \, = \, 
d^{\delta^{x}_{\varepsilon} u}_{\varepsilon}(v,w)$$
\end{prop}

\paragraph{Proof.} 
Use proposition \ref{pplay}. For example, the fact that composition of arrows is
well defined is equivalent with (e) from the mentioned proposition.
Invertibility of arrows is (a), and so on. The isometry claim follows from 
a straight computation. \quad $\square$

This is not a normed groupoid but it can be made into one by adding a 
random walk which embedd the translation groupoid into the transport normed
category.

\subsection{The strong case}

\begin{defi}
A {\bf groupoid strong $\delta$-structure} (or a gs $\delta$-structure) is a
triple  $(G,d, \delta)$ such that $\delta$  is a map assigning to any 
$\varepsilon \in \Gamma$ a transformation $\delta_{\varepsilon}: \, 
dom(\varepsilon) \rightarrow \, im(\varepsilon)$ which satisfies the axioms  
A1, A2 from definition \ref{defgdsweak} and the following axioms 
A3mod and A4: 
\begin{enumerate}
\item[A3mod.] There is a function 
$\displaystyle \bar{d}: U \rightarrow \mathbb{R}$ 
which is the limit  
$$ \lim_{\varepsilon \rightarrow 0} \frac{1}{\mid\varepsilon\mid} \, 
d \, \delta_{\varepsilon} (g)  \, = \, \bar{d}(g)$$
uniformly on bounded sets in the sense of 
definition \ref{defcon}. Moreover, if $\bar{d}(g) = 0$ then $g \in Ob(G)$.  
\end{enumerate}
 \begin{enumerate}
 \item[A4.] the net $\displaystyle \Delta_{\varepsilon}$ converges uniformly on
 bounded sets to a function $\Delta$.  
\end{enumerate}
\label{deltadif}
\end{defi}

\begin{rk}
In the case of a gs $\delta$-structure, notice that A2 and A4 imply that 
the net $\displaystyle dif_{\varepsilon}$ simply converges to $\Delta$,
uniformly on bounded sets.
\label{rema1}
\end{rk}

\begin{prop}
A gs $\delta$-structure is a gw $\delta$-structure. More precisely A1, A2,
A3mod and A4 imply A3 with 
$$\bar{d}(g,h) = \bar{d} \Delta(g,h)$$
\label{pdelta}
\end{prop}

\paragraph{Proof.} 
Indeed, we have: 
$$\frac{1}{\mid\varepsilon\mid} \, d(dif(\delta_{\varepsilon}(g),
\delta_{\varepsilon}(h)) \, = \, \frac{1}{\mid\varepsilon\mid} \, d \,
\delta_{\varepsilon} \, dif_{\varepsilon}(g,h)$$
We reach to the conclusion by using the remark \ref{rema1} and A3mod. 
\quad $\square$

\end{document}